\documentclass[12 pt]{article}
\usepackage{amsmath,amsfonts,amsthm,amssymb,array,mathrsfs,latexsym,paralist, xy, color, stmaryrd}
\usepackage{tikz, soul}
\usepackage{lipsum}

\usepackage{dsfont}

\usetikzlibrary{arrows,automata}

\tikzstyle{V}=[fill=black,circle,scale=0.2, outer sep = 4pt]

\newtheorem{thm}{Theorem}[section]
\newtheorem{prop}[thm]{Proposition}
\newtheorem{cor}[thm]{Corollary}
\newtheorem{lemma}[thm]{Lemma}
\theoremstyle{remark}
\newtheorem{rmk}[thm]{Remark}
\newtheorem{example}[thm]{Example}
\theoremstyle{definition}
\newtheorem{defn}[thm]{Definition}

\newtheorem{assumption}[thm]{Assumption}

\DeclareMathOperator{\Ad}{Ad}

\DeclareMathOperator{\triv}{triv}

\newcommand{\z}{^{(0)}}
\renewcommand{\2}{^{(2)}}
\newcommand{\inv}{^{-1}}

\newcommand{\utimes}{\boxtimes}

\newcommand{\bi}{\begin{itemize}}
\newcommand{\ei}{\end{itemize}}
\newcommand{\be}{\begin{enumerate}}
\newcommand{\ee}{\end{enumerate}}

\newcommand{\C}{\mathbb{C}}
\renewcommand{\S}{\mathcal{S}}

\newcommand{\T}{\mathbb{T}}

\newcommand{\G}{\mathcal{G}}

\newcommand{\R}{\mathbb{R}}
\newcommand{\N}{\mathbb{N}}
\newcommand{\Z}{\mathbb{Z}}

\newcommand{\id}{\operatorname{id}}
\newcommand{\Iso}{\operatorname{Iso}}

\newcommand{\rectimes}{\boxtimes_\theta}
\newcommand{\rtact}{{_\blacktriangleleft}}
\newcommand{\ltact}{{_\blacktriangleright}}
\newcommand{\hact}{{_\diamond}}

\begin{document}

\title{Nontraditional models of $\Gamma$-Cartan pairs}
\author{Jonathan H. Brown and Elizabeth Gillaspy}

\maketitle

\begin{abstract}
This paper explores the tension between multiple models and rigidity for groupoid $C^*$-algebras. We begin by identifying $\Gamma$-Cartan subalgebras $D$ inside twisted groupoid $C^*$-algebras $C^*_r(G, \omega)$,  using similar techniques to those developed in  [DGN$^+$20]. When $D \not= C_0(G\z)$, [BFPR21, Theorem 4.19] then gives another groupoid $H$, and a twist $\Sigma$ over $H$, so that $D \cong C_0(H\z)$ and $C^*_r(G, \omega) \cong C^*_r(H; \Sigma)$.  However, there is a close relationship between  $G$ and $H$.
Theorems \ref{Weyl groupoid} and \ref{thm:Weyl-twist} show how to construct $H$ and $\Sigma$ in terms of $G$ and $\omega$.  We also show, in Theorem \ref{thm:reconstruction-1}, how to reconstruct $G$ from $H$ if we assume the 2-cocycle $\omega$ is trivial.  This latter construction involves a new type of twisting datum (see Proposition \ref{prop:reconstruction-is-groupoid}), which may be of independent interest.
\end{abstract}

\section{Introduction}


Concrete constructions of $C^*$-algebras from other mathematical objects have  led to major advances in both $C^*$-algebra theory and our understanding of the underlying model.  When a given $C^*$-algebra can be described using multiple different models, we obtain a more nuanced and multifaceted picture of the $C^*$-algebra.  
 On the other hand, when a $C^*$-algebra is {\em rigid}, admitting only one  description as $C^*(G)$ for a certain type of objects $G$, the transfer of information between the $C^*$-algebra and the model is more precise. The model gives us a more complete portrait of the $C^*$-algebra, and vice versa.  
 
 Historically, the preponderance of research has focused on rigidity of $C^*$-algebras.  In particular, in the setting of groupoid $C^*$-algebras, Kumjian's groundbreaking work \cite{c*-diagonals} established that, whenever $B$ is a particularly nice maximal Abelian $C^*$-subalgebra of a $C^*$-algebra $A$ (we say $B \subseteq A$ is a {\em $C^*$-diagonal}), then there is a unique principal \'etale groupoid $G,$ and a twist $\Sigma$ over $H$,  so that $A = C^*_r(H; \Sigma)$ is the associated twisted groupoid $C^*$-algebra and $B \cong C_0(H\z)$. Renault \cite{renault-cartan}
subsequently extended Kumjian's work, by replacing ``$C^*$-diagonal'' with the more general notion of ``Cartan subalgebra'' and allowing the groupoids $H$ to be {\em effective}, or topologically principal.   Brownlowe, Carlsen, and Whittaker \cite{brownlowe-carlsen-whittaker} dropped the requirement that the groupoids under consideration be principal or nearly so, at the cost of narrowing their focus to graph $C^*$-algebras.  They showed in \cite[Theorem 5.1]{brownlowe-carlsen-whittaker} that if $E, F$ are directed graphs and $C^*(E) \cong C^*(F)$ via an isomorphism of graph $C^*$-algebras which preserves the canonical diagonal subalgebra,\footnote{Unfortunately, the canonical diagonal subalgebra of a graph $C^*$-algebra is usually not maximal Abelian, and hence is not a $C^*$-diagonal in Kumjian's sense.} then the associated graph groupoids $G_E, G_F$ must be isomorphic.  That is, a graph $C^*$-algebra has a unique graph-groupoid model. A related generalization of Kumjian--Renault theory is found in \cite{BFPR}.  Here the first-named author, together with Fuller, Pitts, and Reznikoff, showed that if a $C^*$-algebra $A$ is graded by a discrete Abelian group $\Gamma$ and $B \subseteq A_0$ is a Cartan subalgebra inside the 0-graded subalgebra $A_0 \subseteq A$, then there is a unique groupoid $H$, with a twist $\Sigma$, so that $A \cong C^*_r(H; \Sigma)$ and $B \cong C_0(H\z)$.

On the other hand, Latr\'emoli\`ere and Packer's analysis of the noncommutative solenoids \cite{LP} is
 an example of a situation where identifying multiple group(oid) models for a  $C^*$-algebra led to key insights.  The noncommutative solenoids can be described as a twisted group $C^*$-algebra or as the $C^*$-algebra of a transformation groupoid.  Latr\'emoli\`ere and Packer used the first description to compute the $K$-theory of the noncommutative solenoids, and the second one to characterize simplicity of these $C^*$-algebras. 

In this paper, we investigate the tension between multiple models and rigidity for groupoid $C^*$-algebras, building on the aforementioned \cite{BFPR} and the prior work of the second-named author with Duwenig, Norton, Rezni-koff, and Wright \cite{DGNRW, DGN}. Theorem 3.1 of \cite{DGNRW} identifies situations when, given  a non-effective groupoid $G$, we can still find a Cartan subalgebra inside the twisted groupoid $C^*$-algebra $C^*_r(G, \omega)$.  Renault's theorem \cite{renault-cartan} then implies that $C^*_r(G, \omega) \cong C^*_r(H; \Sigma)$ for some effective groupoid $H \not= G$ and a twist $\Sigma$ over $H$; that is, we have multiple groupoid models for the same $C^*$-algebra.  The authors of \cite{DGN}  show how to construct $H$ from $G$, and identify the twist $\Sigma$ in certain cases.

Here,  we identify in Theorem \ref{gamma cartan exists} sufficient conditions for a twisted groupoid $C^*$-algebra $C^*_r(G, \omega)$ to admit a $\Gamma$-Cartan subalgebra in the sense of \cite{BFPR}. 
As in \cite{DGNRW}, the $\Gamma$-Cartan subalgebras that we identify arise from the dynamics of $G$, so we will refer to them as dynamical $\Gamma$-Cartan subalgebras. 
These dynamical $\Gamma$-Cartan subalgebra will generally be strictly larger than $C_0(G\z)$, so $G$ cannot be the groupoid $H$ constructed by Brown et al. 
As in \cite{DGN}, one major focus of this paper is to understand the relationship between the original data $(G, \omega)$ and the groupoid-twist pair $(H, \Sigma)$ guaranteed by \cite{BFPR} (see Theorems \ref{Weyl groupoid} and \ref{thm:Weyl-twist} below, as well as Section \ref{sec:reconstruction}).  

The Weyl groupoid $H$ which  \cite{BFPR} associate to a $\Gamma$-Cartan subalgebra is a relatively rigid object; for one thing, $H$ must admit a homomorphism $c: H \to \Gamma$ with $c^{-1}(0)$ effective. 
Indeed, we have rigidity in the sense of the first paragraph.  By \cite[Corollary 6.3]{BFPR}, a $C^*$-algebra $A$ admitting a $\Gamma$-Cartan subalgebra can be described by exactly one pair $(H, c)$ with $c^{-1}(0)$ effective.  

Despite this rigidity, hindsight enables us to identify examples that have already appeared in the literature of groupoid $C^*$-algebras with ``unexpected'' $\Gamma$-Cartan $C^*$-algebras. 
If $G$ is a group and $\omega$ is a $2$-cocycle, and $c: G \to \Gamma$ is a homomorphism,  $c\inv(0)$ is effective if and only if $c$ is injective.  It follows from 
\cite{BFPR} that, for a group $G$, $\C = C_0(G\z)$ is $\Gamma$-Cartan in $C^*_r(G, \omega)$ exactly when   $\Gamma=G$ and $c=\id$. 
However, Theorem \ref{gamma cartan exists} below 
tells us that, for example, both the rotation algebras $C^*_r(\Z^2, \omega_\theta)$ and the noncommutative solenoids 
 $C^*_r(G, \omega)$  of
\cite{latremoliere-packer, LP} admit   a $\Gamma$-Cartan subalgebra for $\Gamma \not= G$. 
We discuss this in more detail in Examples \ref{ex:nc-tori} and 
\ref{ex: solenoid} below.


Having identified, in Theorem \ref{gamma cartan exists}, twisted groupoid $C^*$-algebras $C^*_r(G, \omega)$ which admit unexpected or  ``nontraditional'' $\Gamma$-Cartan subalgebras,   our next goal is to identify which Weyl groupoids $H$  arise in this fashion.  From Theorem \ref{Weyl groupoid},  we know that $H = G/S \rtimes \widehat S$ for a subgroupoid $S$ of $G$.  Abstracting this relationship, Assumption \ref{assumption-on-actions} (see also Proposition \ref{prop:weyl-satisfies}) lists structural properties that the Weyl groupoids arising from Theorem \ref{gamma cartan exists} must possess.  Indeed, if $H$ satisfies Assumption \ref{assumption-on-actions} and we assume the original 2-cocycle $\omega$ was trivial, Theorem \ref{thm:reconstruction-1} shows that we can reconstruct the original $G$ from $H$, 
together with an additional twisting datum $\theta$, which has something of the flavor of a 2-cocycle but is (to our knowledge) a novel construction.  We hope that this new method of constructing one groupoid from another  will be of independent interest.

This paper is organized as follows.  We review the basics of (twisted) groupoid $C^*$-algebras and $\Gamma$-Cartan subalgebras in Section \ref{sec:prelim}. Section \ref{sec:main-results} contains our first main theoretical results.  We describe in Theorem \ref{gamma cartan exists} when a groupoid $C^*$-algebra admits a dynamical $\Gamma$-Cartan subalgebra, and we describe the associated Weyl groupoid and (under additional hypotheses) Weyl twist in Theorems \ref{Weyl groupoid} and \ref{thm:Weyl-twist} respectively.  Section \ref{sec:applications} applies these results to the setting of twisted group $C^*$-algebras, in particular to the noncommutative solenoids of \cite{LP}.  Finally, Section \ref{sec:reconstruction} begins by identifying (see Assumption \ref{assumption-on-actions} and Proposition \ref{prop:weyl-satisfies})  which  groupoids could possibly arise as the Weyl groupoid of a dynamical $\Gamma$-Cartan pair, and closes by showing  (in Theorems \ref{thm:uniqueness} and  \ref{thm:reconstruction-1}) that in many cases, we can indeed recover the original groupoid of the dynamical $\Gamma$-pair from these Weyl groupoids.

\subsection*{Acknowledgments}
This research was partially supported by the National Science Foundation (grant 1800479 to E.G.).  JHB thanks EG for hospitality during visits to the University of Montana.

\section{Preliminaries}
\label{sec:prelim}

\begin{defn} \cite{renault-cartan}
A {\em Cartan subalgebra} of a $C^*$-algebra $A$ is a maximal Abelian subalgebra $B$ such that:
\begin{enumerate}
\item There is a faithful conditional expectation $E: A \to B$;
\item The normalizers of $B$ in $A$, $\{ n \in A: n ^* b n , n b n^* \in B \text{ for all } b \in B\}$, generate $A$ as a $C^*$-algebra;
\item $B$ contains a countable approximate unit for $A$.
\end{enumerate}
\end{defn}
Renault proved that if $B \subseteq A$ is Cartan, then $A = C^*_r(G, \Sigma)$ and $B = C_0(G\z)$ for a unique effective groupoid $G$ and a twist $\Sigma$ over $G$ (see Section \ref{sec:gpoid} for definitions).

\begin{defn}\cite[Definition 2.1]{BFPR}
	A $C^*$-algebra $A$ is {\em topologically graded} by a discrete Abelian group $\Gamma$ if there are linearly independent closed linear subspaces $A_s$ of $A$, for $s \in \Gamma$, such that:
	\begin{enumerate}
		\item $A_s A_t \subseteq A_{s+t}$  for all $s, t\in \Gamma$;
		\item $A_t^* = A_{-t}$;
		\item $A$ is densely spanned by $\{A_s\}_{s\in \Gamma}$; and
		\item there is a faithful conditional expectation from $A$ to $A_0$.
	\end{enumerate} 
\label{def:grading}
\end{defn}

\begin{defn}\cite[Definition 3.9]{BFPR}
Let $A$ be a  $C^*$-algebra topologically graded by a discrete Abelian group $\Gamma$ and $D$
	an Abelian $C^*$-subalgebra of $A_0$. We say the pair $(A,D)$ is {\em $\Gamma$-Cartan} if $D$ is Cartan in $A_0$ and $ N(A,D)$ spans a dense subset of $A$.
\end{defn}

Brown et al.~proved in \cite[Theorem 4.36]{BFPR} that if $D \subseteq A$ is a $\Gamma$-Cartan subalgebra, then there is a unique triple $(G, \Sigma, c)$, consisting of an \'etale groupoid $G$ [not necessarily  effective]; a twist $\Sigma$ over $G$; and a 1-cocycle or group(oid) homomorphism $c: \Sigma \to \Gamma$, such that $C^*_r(G; \Sigma) \cong A$; $C_0(G\z) \cong D$; and the $\Gamma$-grading comes from $c$, in the sense that $A_s$ consists of functions supported on $c^{-1}(s) \subseteq \Sigma$.

\subsection{Groupoids}
\label{sec:gpoid}
Recall that a {\em groupoid} is a small category $G$ in which every element $\gamma \in G$ has an inverse $\gamma^{-1}$, in the sense that $(\gamma^{-1})^{-1} = \gamma$,  $\gamma \gamma^{-1} \gamma = \gamma$ and $\gamma^{-1} \gamma  \gamma^{-1} = \gamma^{-1}$.  We write $G\z = \{ \gamma \gamma ^{-1}: \gamma \in G\}$ for the space of {\em units} of $G$, and define $r(\gamma) = \gamma \gamma^{-1}$ and $s(\gamma) = \gamma^{-1} \gamma$.  The set of {\em composable pairs} is $G\2 = \{ (\gamma,\eta)\in G: s(\gamma ) = r(\eta  )\}$; if $(\gamma,\eta)\in  G\2$ then $\gamma \eta  \in G$.  Moreover, this multiplication is associative.

For $u \in G\z$, we write $G^u = \{ \gamma \in G: r(\gamma ) = u\}$ and $G_u = \{ \gamma \in G: s(\gamma ) =u\}$.
For a groupoid  $G$, we write $\Iso(G) = \{ \gamma  \in G: r(\gamma ) = s(\gamma )\}$.

\begin{defn}
	A subgroupoid $S \leq G$ is {\em normal} if, for all $\gamma \in G$, $\gamma^{-1} S \gamma \subseteq S$.
\end{defn}

We will generally assume that our groupoids are {\em topological} groupoids, in the sense that $G$ has a (locally compact Hausdorff) topology with respect to which the multiplication and inversion maps are continuous.  In fact, many of the groupoids in this paper will be {\em \'etale}, meaning that $r$ will be a local homeomorphism.  
Equivalently, an \'etale groupoid $G$ always has a basis consisting of open bisections, where a {\em bisection} is a subset $B\subseteq G$ such that there exists an open set $U \supseteq B$ for which $r|_U, s|_U$ are both homeomorphisms. 
We say that $G$ is {\em effective} if $\Iso(G)^\circ = G\z$.

\begin{defn}
A {\em twist} $\Sigma$ over a groupoid $G$ is a groupoid $\Sigma$ with unit space $G\z$ such that $\Sigma$ is a principal $\T$-bundle over $G$, and the  groupoid multiplication commutes with the $\T$-action:  if $(\sigma, \tau) \in \Sigma\2$, then for any $z, w \in \T$ we have $(\sigma z, \tau w) \in \Sigma\2$ and $\sigma z \cdot \tau w = (\sigma \tau) zw$.  The most accessible examples of twists come from 2-cocycles.

A {\em (continuous) 2-cocycle} $\omega $ on a topological groupoid $G$ is a continuous function $\omega: G\2 \to \T$ satisfying the {\em cocycle condition:} whenever $(g,h), (h, k) \in G\2$,
\[ \omega(g, hk) \omega (h, k) = \omega (gh, k) \omega (g,h).\]
Given a  2-cocycle $\omega$ on $G$, we define a twist $\Sigma_\omega$ over $G$ by setting $\Sigma_\omega = G \times \T$ as topological spaces, and defining the groupoid structure on $\Sigma$ by 
\[ (g, z)(h, w) = (gh, \omega(g,h) zw); \quad (g,z)^{-1} = (g^{-1}, \overline{z \omega(g, g^{-1})}).\] 
\end{defn}

Given a locally compact Hausdorff \'etale groupoid $G$ and a 2-cocycle $\omega $ on $G$, we can define a twisted $*$-algebra structure on $C_c(G)$:  for $\phi, \psi \in C_c(G)$, 
\begin{equation}
\phi * \psi (g) := \sum_{(h,k) \in G\2: hk = g} \phi(h) \psi(k) \omega(h,k); \qquad \psi^*(g) = \overline{\psi(g^{-1}) \omega(g, g^{-1})}.
\label{eq:conv-mult}
\end{equation}
For each $u \in G\z$, we have a left-regular representation $\pi_u$ of this convolution algebra $C_c(G, \omega)$ on $\ell^2(G_u)$: 
\[ \pi_u(\phi) (\xi) := \phi * \xi.\]
One can compute that $\pi_u(\psi)^* = \pi_u(\psi^*)$.

 The {\em reduced groupoid $C^*$-algebra} $C^*_r(G,\omega)$ is the completion of  $C_c(G,\omega)$ in the norm 
 \[ \| \psi\|_r := \sup_{u \in G\z} \| \pi_u(\psi)\|.\]
 
 One obtains a similar definition for $C^*_r(G, \omega)$ if $G$ is not \'etale, using integration against a Haar measure  on $G$ instead of the sum in Equation \eqref{eq:conv-mult}; cf.~\cite{williams-gpoid}.  One can also define $C^*_r(G; \Sigma)$ if $\Sigma$ is a general twist over $G$, in such a way that $C^*_r(G; \Sigma_\omega) \cong C^*_r(G, \omega)$; 
see \cite[Remark 5.1.6]{sims-notes}.


Given a locally compact Hausdorff  groupoid $G$ (with unit space $G\z$) and a topological space $Z$ equipped with a continuous open surjection $p: Z \to G\z$, we can form the fibre products $G \mathbin{_s * _p} 
 Z$ and $Z\mathbin{_p * _r} G$.  We say that $Z$ is a $G$-space and $p$ is the moment map. If we have a continuous map $\alpha: G \mathbin{_s * _p} Z \to Z$ which satisfies $p(\alpha(g,z)) = r(g); \alpha(gh, z) = \alpha(g, \alpha(h,z));$ and  $\alpha(u, z) = z$ for all $u \in G\z$; then we say $G$ acts on the left of $Z$.
Similarly, a right action of $G$ on $Z$ is given by a continuous map $\beta: Z {}_p *_r G \to Z$ satisfying $p(\beta(z,g)) = s(g); \beta(z, gh) = \beta(\beta(z, g), h) );$ and $\beta(z, u) = z$ for all $u \in G\z$.  We will sometimes write $ \alpha_g(z) $ or $ g.z$ instead of $\alpha(g,z)$.

A left action is free if $\alpha(g, z) = z$ implies $g = p(z)$; it's proper if the map $G\mathbin{_s * _p} Z \to Z \times Z$ given by $(g, z) \mapsto (\alpha(g,z), z)$ is proper in the sense that preimages of compact sets are compact.   Equivalently \cite[Proposition 2.17]{williams-gpoid}, a left action is proper if and only if whenever we have nets $\{z_i\}_{i\in I} \subseteq Z$ and $\{g_i\}_{i\in I} \subseteq  G$ such that $z_i \to z$ and $g_i . z_i \to w$,  then $\{ g_i\}_{i\in I}$ has a convergent subnet.

It follows in particular from
\cite[Proposition~2.18]{williams-gpoid}  that if $G$ is a locally compact Hausdorff \'etale groupoid and $X$ is a proper (left) $G$-space, then the quotient space $G\backslash X$ is also locally compact and Hausdorff.

If a subgroupoid $S \leq G$ is {\em wide}, i.e., $S\z  = G\z$, then $S$ acts on the left of $G$ (taking the moment map to be $p = r$) and on the right of $G$ (taking the moment map to be $p = s$).

{
\begin{prop}
\label{prop:proper-action}
Let $G$ be a locally compact Hausdorff groupoid and let $S\leq G$ be a closed wide subgroupoid.  Then $G$ is a proper (right)  $S$-space.
\end{prop}
\begin{proof}
Suppose that $(\gamma_i)_i \to \gamma \in G$ and that $(\rho_i) \subseteq S$ satisfies $(\gamma_i \cdot \rho_i)_i \to \eta \in G.$  Since $S$ is closed, and  $\rho_i = \gamma_i^{-1} \gamma_i \rho_i \to \gamma^{-1} \eta$, we have $\lim_i \rho_i = \gamma^{-1}\eta \in S$.  We conclude that $G$ is a proper $S$-space.
\end{proof}
}

\section{Groupoids with $\Gamma$-Cartan subalgebras}
\label{sec:main-results}

In this section, we identify conditions under which a subgroupoid $S$ of $G$ will give rise to a $\Gamma$-Cartan subalgebra of the twisted groupoid $C^*$-algebra $C^*_r(G, \omega)$.  We also give an explicit description of the associated Weyl groupoid (Theorem \ref{Weyl groupoid}) and, under mild hypotheses, of the Weyl twist (Theorem \ref{thm:Weyl-twist}).
The arguments in this section closely follow those of \cite{DGN}.

We begin with a key hypothesis; cf.~\cite[p.~8]{DGNRW} for more discussion of this property.

\begin{defn}
\label{def:imm-cent}
Let $S \leq T$ be group bundles.  We say $S$ is {\em immediately centralizing} (in $T$) if, for all $t \in T$ and all integers $k \geq 1$, the statement 
\[ (\forall \, s \in S_{p(t)})(\exists \, 1 \leq n \leq k)(t s^n = s^n t)\]
implies that $ts = st$ for all $s \in S_{p(t)}$.
\end{defn}

Intuitively, $S$ is immediately centralizing in $T$ if whenever an element $t \in T$  commutes with all sufficiently large powers of all $s \in S_{p(t)}$, then $t$ commutes with $S_{p(t)}$.

\begin{thm}\label{gamma cartan exists}
Suppose $G$ is a Hausdorff \'etale groupoid, $\Gamma$ is a discrete Abelian group, $c: G\to \Gamma$ is a homomorphism and $\omega$ is a continuous $2$-cocycle on $G$. Suppose further that $S\subseteq \Iso(c\inv(0))$ is maximal among all open Abelian subgroupoids of $\Iso(c\inv(0))$ on which $\omega$ is symmetric.  If $S$ is closed and normal in $G$ and immediately centralizing in $\Iso(c^{-1}(0))$, then $C^*_r(S,\omega)\subseteq C^*_r(G,\omega)$ is $\Gamma$-Cartan.
\end{thm}

\begin{rmk}
While this may look like a long list of hypotheses, they are nearly all necessary in order for a subgroupoid $S\leq G$ to give rise to a $\Gamma$-Cartan subalgebra of $C^*_r(G, \omega)$: 
\begin{itemize}
\item (cf.~ \cite[Lemma 2.19]{BFPR}) We need $S \leq G$ to be open in order to obtain the inclusion  $C^*_r(S, \omega) \subseteq C^*_r(G, \omega)$.  Note that if $S$ is not open, $C_c(S)$ does not naturally embed in $C_c(G)$. 
\item Having $S$ be Abelian and $\omega|_{S \times S}$ be symmetric are necessary in order for $C^*_r(S, \omega)$ to be Abelian. (This is easily checked by choosing two elements $a, b \in S$ with $ab \not= ba$ and/or $\omega(a,b) \not= \omega(b,a)$, and multiplying two bump functions $f_a, f_b \in C_c(S, \omega) \subseteq C^*_r(S, \omega)$ supported in [disjoint] open-bisection neighborhoods $U_a, U_b$ of $a,b$ respectively.)
\item We need maximality of $S$ among open Abelian subgroupoids of $\Iso(c^{-1}(0))$ on which $\omega$ is symmetric in order for $C^*_r(S, \omega)$ to be maximal Abelian in $C^*_r(G, \omega)$. (If $S$ were not maximal, there would be a strictly larger subgroupoid $T \leq \Iso(c^{-1}(0)$, so that $C_c(S, \omega) \subsetneqq C_c(T, \omega) \subseteq C^*_r(T, \omega)$, where this larger algebra is also Abelian (cf.~\cite[Lemma 3.5]{DGNRW}).) 
\item (cf.~\cite[Lemma~3.4]{BEFPR}) The need for a conditional expectation of $C^*_r(G, \omega)$ onto $C^*_r(S, \omega)$ forces $S$ to be closed. 
\item Normality of $S$ is needed in order for normalizers of $C^*_r(S, \omega)$ to generate $C^*_r(G, \omega)$. (If $a \in S, \gamma \in G$ but $\gamma a \gamma^{-1}\not\in S$, then for all sufficiently small open-bisection neighborhoods $U_\gamma$ of $\gamma$,  a bump function $f_\gamma$ supported in $U_\gamma$ will not normalize $C^*_r(S)$.)
\end{itemize}
We note, however, that a  failure to find a subgroupoid $S \leq G$ satisfying the hypotheses of Theorem \ref{gamma cartan exists} does not preclude the existence of Cartan subalgebras of $C^*_r(G, \omega)$ which arise by  non-dynamical  means.
\end{rmk}

One important step in the proof of Theorem \ref{gamma cartan exists}
is the following:
\begin{prop}
\label{prop:revised-dgnrw}
Suppose $H$ is a Hausdorff \'etale groupoid, $\omega$ is a continuous 2-cocycle on $H$, and $S \subseteq \Iso(H)$ is maximal among open Abelian subgroupoids  of $\Iso(H)$ on which $\omega$ is symmetric.  If $S$ is closed and normal in $H$ and immediately centralizing in $\Iso(H)$, then $C^*_r(S, \omega)$ is Cartan in $C^*_r(H, \omega)$.
\end{prop}
\begin{proof}
The difference between this proposition and \cite[Theorem 3.1]{DGNRW} is the location of the adjective ``open.'' In order to apply the proof of that theorem in our setting, we need to establish the following analogue of \cite[Lemma 3.6]{DGNRW}:

\begin{center}  ``Suppose $S$ is maximal among Abelian open subgroupoids of $\Iso(H)$ on which $\omega$ is symmetric. If $U \subseteq \Iso(H)$ is open and whenever $(\eta, s) \in ( U \times S)\cap H\2$ we have $\eta s = s \eta$ and $\omega(s,\eta) = \omega(\eta, s)$, then $U \subseteq S.$'' 
\end{center}

To that end, we begin by observing that since $H$ is \'etale,  such an open set $U$ must contain a bisection $W_\eta$ around each element $\eta \in U$.  As $W_\eta$ is a bisection and $S \subseteq \Iso(H)$,  elements in the groupoid generated by $S$ and $W_\eta$ are all of the form $\rho  \gamma^k$ for some $\rho \in S,  \gamma \in W_\eta$. (There is no need to worry about products in other orders, as we assumed that $W_\eta\subseteq U$ and  that $U, S$ commute.)  Furthermore, if $(\rho  \gamma^k, t \xi^n) \in H\2$, the fact that $S$ is a subgroupoid of $\Iso(H)$ implies that $s(\gamma) = r(t) = s(t) = r(\xi)$ and hence the fact that $W_\eta$ is a bisection implies that $\xi = \gamma$.  Now, the inductive argument of \cite[Lemma 3.6]{DGNRW} implies that $\omega$ is symmetric on the groupoid $K \subseteq \Iso(H)$ generated by $S$ and $W_\eta$.  Moreover, this groupoid is open. The maximality of $S$ therefore implies that $S \supseteq K$, ie, $W_\eta \subseteq S$ for all $\eta$.  Hence, $U \subseteq S$.

Applying this statement instead of \cite[Lemma 3.6]{DGNRW} wherever necessary in the proof of \cite[Theorem 3.1]{DGNRW}, we obtain the desired conclusion.
\end{proof}

\begin{proof}[Proof of Theorem \ref{gamma cartan exists}]
We view $A := C^*_r(G,\omega)$ as a $\Gamma$-graded $C^*$-algebra via $c$: Define 
\[ A_t = {  \{ f \in C_c(G, \omega): f(\gamma) \not= 0 \implies c(\gamma) = t\} } = { C_c(G, \omega)|_{ c^{-1}(t)}}.\]
 The fact that $c$ is continuous implies that for any $f\in C_c(G, \omega), \ f|_{ c^{-1}(t)} \in C_c(G, \omega)$. 
Observe that whenever $f \in C_c(G, \omega)|_{ c^{-1}(t)}, g \in C_c(G, \omega)|_{ c^{-1}(s)}$ and $f * g(\gamma) \not= 0$, 
\begin{align*}
	f*g(\gamma) &= \sum_{r(\eta) = s(\gamma)} f(\gamma \eta) g(\eta^{-1}) \omega(\gamma \eta, \eta^{-1})\\
	&= \sum_{\eta: r(\eta) = s(\gamma), c(\eta) = -s } f(\gamma \eta) g(\eta^{-1}) \omega(\gamma \eta, \eta^{-1})
\end{align*}
is nonzero iff $c(\gamma) = t +s$; that is, $A_t A_s \subseteq A_{t+s}$. Properties (2) and (3) of Definition \ref{def:grading} are also easily verified.  The proof that $\Phi_0: C^*_r(G, \omega) \to A_0$, given on $C_c(G, \omega)$ by restriction, is a faithful conditional expectation proceeds exactly as in \cite[Proposition 3.13]{DGNRW}, using the observation that $c(u) = 0$ for any unit $u\in G\z$. 

Proposition \ref{prop:revised-dgnrw} implies that $C^*_r(S, \omega)$ is Cartan in $A_0 = C^*_r(c^{-1}(0), \omega)$, completing the proof.
\end{proof}

Given a $\Gamma$-Cartan subalgebra $D \subseteq A$, we have  a faithful conditional expectation $\Delta: A \to D$, given by composing the faithful conditional expectation $A_0 \to D$ with the faithful conditional expectation $A \to A_0$.  

\begin{prop}
If a $\Gamma$-Cartan subalgebra $D$ of $A = C^*_r(G, \omega)$ satisfies $D = C^*_r(S,\omega) \cong C_0(\widehat{C^*_r(S,\omega)})$  as in Theorem \ref{gamma cartan exists}, then the conditional expectation $\Delta: A \to D$ is given on $C_c(G, \omega) \subseteq C^*_r(G, \omega)$ by 
\[( \Delta f)x = \sum_{a \in S_{p(x)} }x(a) f(a).\]
Here $x \in \widehat{C^*_r(S,\omega)}$, which is a bundle of groups fibred over $S^{(0)}$; we write $p: \widehat{C^*_r(S,\omega)} \to S^{(0)}$ for the projection map.
\end{prop}

\begin{proof}
As we noted above, the proof of \cite[Proposition 3.13]{DGNRW} establishes that the map $\Phi_0: A \to A_0$ which is given on $C_c(G, \omega)$ by restriction to $c^{-1}(0) \subseteq G$ is a faithful conditional expectation. 
Therefore, the formula given in \cite[Lemma 3.2]{DGN} 
for the faithful conditional expectation associated with the Cartan inclusion $C^*_r(S, \omega) \subseteq A_0 = C^*_r(c^{-1}(0), \omega)$, which is valid for functions $f \in C_c(c^{-1}(0), \omega)$, also yields the formula for $\Delta$ when applied to  functions $f \in C_c(G, \omega)$. 
\end{proof}

{Recall from \cite{BFPR} that a $\Gamma$-Cartan subalgebra $D$ in a $C^*$-algebra $A$ gives rise to a groupoid description of $A$: we have $A \cong C^*_r(G_W; \Sigma_W)$ and $D \cong C_0(G_W\z)$ where $G_W$ is the Weyl groupoid and $\Sigma_W$ is the Weyl twist of the $\Gamma$-Cartan pair.  Consequently, $G_W\z = \widehat D$.
\begin{defn}
\label{def:weyl-gpoid} \cite[Definition 4.9]{BFPR}

The Weyl groupoid $G_W$ and Weyl twist $\Sigma_W$ of the dynamical $\Gamma$-Cartan pair arising from $S \leq G$ are 
	\begin{align*}
	G_W & = \{ [n, x]: x \in \widehat{C^*_r(S, \omega)}, n \in N(C^*_r(S, \omega))\}, \\
	  \Sigma_W& = \{ \llbracket n, x \rrbracket : x \in \widehat{C^*_r(S, \omega)}, n \in N(C^*_r(S, \omega))\},
	  \end{align*}
	where $[n, x] = [m, y]$ iff $x = y$ and $\Delta(m^*n)(x) \not= 0$, and $\llbracket n, x \rrbracket = \llbracket m, y\rrbracket $ iff $x = y $ and $\Delta(m^*n)(x) > 0$.  Moreover (cf.~\cite[Lemma 4.28]{BFPR}), for $n \in C_c(G, \omega)$ supported in a bisection   and $\mathcal O\subseteq \T$ open, the sets 
	\[ Z(n, \mathcal O) := \{ \llbracket \lambda n, x \rrbracket: \Delta(n^* n)(x) > 0, \lambda \in \mathcal O \} \]
	form a basis for the
	 topology on $\Sigma_W$.
	
	If $f \in D  = C_0( \widehat{C^*_r(S, \omega)})$, we write $\text{supp}\, f$ for the {\em open} support of $f$: that is,
	\begin{equation} 
	\text{supp}\, f  = \{ x \in \widehat{C^*_r(S, \omega))}: f(x) \not= 0 \}.
	\label{eq:support-defn}
	\end{equation} 
	 For each $n \in N(C^*_r(S, \omega))$, the proof of \cite[Proposition 1.6]{c*-diagonals} establishes the existence of a unique partially defined action $\alpha_n: \text{supp}\, \Delta(n^*n) \to \text{supp}\, \Delta(nn^*)$, such that for each $f \in D$ and each $x \in \text{supp}\, \Delta(n^*n)$, 
	 \[n^*fn(x) = f(\alpha_n(x)) n^*n(x).\] 
	 Among other things, this action allows us to define the multiplication in $G_W$ and $\Sigma_W$:
	the product $[n, x][m, x']$ is defined iff $x = \alpha_m(x')$,  in which case 
	\[ [n, \alpha_m(x)][m, x] = [nm, x].\]  Thus, $s([m,x]) = x$.
\end{defn}
}

Theorem \ref{Weyl groupoid} below establishes an alternative description of the Weyl groupoid $G_W$ associated to the $\Gamma$-Cartan pairs of Theorem \ref{gamma cartan exists}.  To be precise,  $G_W$ is homeomorphic to the fibre product  $G/S *_p \widehat S$.  (The quotient groupoid $G/S$ is again locally compact Hausdorff, thanks to Proposition \ref{prop:proper-action}; therefore, so is the fibre product.)  Moreover, the proof of \cite[Proposition 4.1]{DGN} goes through without change, to show that 
\begin{align} 
	\label{eq:Weyl-gpoid-action}
	[\gamma] . x(a) &= \overline{\omega(\gamma, \gamma^{-1})} \omega(\gamma^{-1}, a) \omega(\gamma^{-1}a, \gamma) x(\gamma^{-1} a \gamma)\\
	&  = \omega(\gamma^{-1}, a \gamma) \overline{\omega(a\gamma, \gamma^{-1})} x(\gamma^{-1}a \gamma) \notag
\end{align}
describes a well-defined action of $G/S$ on $\widehat S$.  We write $G/S \rtimes \widehat S$ for the associated action groupoid.  That is, 
\[ G/S\rtimes\widehat S = \{ ([\gamma], x): p(x) = s(\gamma)\},
\]
 where $r([\gamma], x)=[\gamma] . x$, $s([\gamma], x)=x$, multiplication is given by $([\gamma], [\eta].x)([\eta],x)=([\gamma\eta],x)$, and inversion is $([\gamma], x)\inv=([\gamma\inv], [\gamma].x)$.

\begin{thm}\label{Weyl groupoid}Suppose 
that $G, \Gamma, c, \omega, S$ satisfy the hypotheses of Theorem \ref{gamma cartan exists}. Then the Weyl groupoid $G_W$ associated to the $\Gamma$-Cartan pair $(C^*_r(G, \omega), C^*_r(S, \omega))$ is isomorphic to $G/S \rtimes \widehat S$. 

The isomorphism $\Psi: G/S \rtimes \widehat S \to G_W$  is  given by $\Psi([\gamma], x) := [n, x]$, where $n \in N_h \subseteq A = \overline{C_c(G)}$ satisfies $n(\gamma) \not= 0$. 
\end{thm}
\begin{proof}
The proof of this Theorem follows the proof of \cite[Theorem 4.6]{DGN}. 
As $S$ is a bundle of Abelian groups on which $\omega$ is symmetric, \cite[Section 3]{DGN} describes how to view $\widehat S$ as the Gelfand dual of the Abelian $C^*$-algebra $ D = C^*_r(S, \omega)$.  
For \cite[Lemma 4.2]{DGNRW}, a straightforward computation establishes  that $\{f \in C_c(G, \omega): f \text{ supported in a bisection}\} \subseteq N(A, D)$, and we know from  \cite{errata} that \cite[Proposition 4.1]{DGNRW} holds for the Weyl groupoids of $\Gamma$-Cartan pairs as well.  Similarly, \cite[Lemma 4.3, Proposition 4.4, Proposition 4.5]{DGN} hold without change, as the proofs of these results do not rely on $S$ being maximal Abelian in $G$.
\end{proof}

Observe that for $([p(x)], x) \in (G /S \rtimes \widehat S)\z$, we have $\Psi([p(x)], x) = [f, x]$ for any function $f \in C_0(G\z)$ which is nonzero near $p(x)$.  In particular,  on $\widehat S = \widehat{C^*_r(S, \omega)}$, $\Psi= \text{id}$.

In certain situations, we can also describe the Weyl twist $\Sigma_W$ associated to $(C^*_r(G, \omega), C^*_r(S, \omega))$.

\begin{thm}
\label{thm:Weyl-twist}
Suppose 
that $G, \Gamma, c, \omega, S$ are as in the previous Theorem.  
Suppose further that $\frak s : G/S\to G$ is a continuous section.  Then the associated Weyl twist is topologically trivial. In particular, $C^*_r(G, \omega ) \cong C^*_r(G/S \rtimes \widehat S, C)$ where the continuous 2-cocycle  $C$ is given by the formula 
\begin{align*}
C(([\gamma], [\eta].x),([\eta], x)) &= \\ 
= x(\frak s([\gamma \eta])^{-1}\frak s([\gamma]) & \frak s([\eta]) ) \omega(\frak s([\gamma\eta])^{-1}, \frak s([\gamma]) \frak s([\eta])) \times \\
& \qquad \times  \overline{\omega(\frak s([\gamma \eta])^{-1}, \frak s([\gamma \eta]))} \omega(\frak s([\gamma]), \frak s([\eta])) \\
 =  x(\frak s([\gamma \eta])^{-1}\frak s([\gamma]) & \frak s([\eta]) ) \overline{\omega(\frak s([\gamma \eta]), \frak s([\gamma \eta])^{-1} \frak s([\gamma]) \frak s([ \eta]))} \omega(\frak s([\gamma], \frak s([\eta])).
 \end{align*}
 \end{thm}

\begin{proof}
	Recall the definition of the Weyl groupoid and Weyl twist from Definition \ref{def:weyl-gpoid}.	
We will first define a section $\psi_{\frak s}: G_W \to \Sigma_W$, which we will then use to construct $C$.   For each $[n,x] \in G_W$,  choose a function $f_{n,x} \in C_c(G, \omega)$, supported in a bisection, such that $f_{n,x}(\frak s([\gamma])) > 0$, where $[n, x] \leftrightarrow ([\gamma], x)$ via the identification of Theorem \ref{Weyl groupoid}. Then the map $\psi_{\frak s}([n,x]) = \llbracket f_{n,x}, x \rrbracket$ defines a section of the Weyl twist $\Sigma_W \to G_W$.  

To see that $\psi_{\frak s}$ is well defined, note that if $f \in C_c(G, \omega )$ is also supported in a bisection and  $f(\frak s([\gamma])) > 0$, then $\frak s([\gamma])$ is the unique element $\eta$ of $G_{p(x)}$ with $f(\eta) \not= 0$.  Consequently, 
\begin{align*}
	\Delta( f^* f_{n,x}) (x) & = \sum_{a \in S_{p(x)} } x(a) f^* f_{n,x}(a)\\
	&  = \sum_{a \in S_{p(x)}}   \sum_{s(\eta) =p( x)} x(a) \overline{f(\eta) \omega(\eta^{-1}, \eta)} f_{n,x}(\eta a) \omega(\eta^{-1}, \eta a) \\
	& = \overline{f(\frak s([\gamma]))} f_{n,x}(\frak s([\gamma]))  > 0,
\end{align*}
since in order to have both $f(\eta) $ and $f_{n,x}(\eta a)$ nonzero for $a \in S_{p(x)}$, we must have $a = p(x) \in G\z$.
It follows that  $\llbracket f, x \rrbracket = \llbracket f_{n,x}, x \rrbracket $.

To see that $\psi_{\frak s}$ is continuous,  suppose that $[n_i, x_i ] \to [n,x]$.  Under the isomorphism $G/S \rtimes \widehat S \cong G_W$, we can equivalently describe these groupoid elements as  $([\gamma_i], x_i) \to ([\gamma], x)$.  Since $\frak s$ is continuous,   $\frak s([\gamma_i]) \to \frak s([\gamma])$, so there exists $I$ such that $i \geq I$ implies $\frak s([\gamma_i]) \in \text{supp}(f_{n,x})$.  
It follows that, for $i \geq I$, defining 
\[ \lambda_i := \frac{ \overline{f_{n, x}(\frak s[\gamma_i])}}{| f_{n,x}(\frak s[\gamma_i])|}\]
yields $\llbracket \lambda_i f_{n,x}, x \rrbracket = \llbracket f_{n_i, x_i}, x \rrbracket$: as everything in sight is supported in a bisection, a similar computation to the one in the second paragraph of the proof establishes that, since $\frak s[\gamma_i] \in \text{supp}(f_{n,x})$ and $s(\frak s[\gamma_i] a ) = s(\frak s[\gamma_i])$ for any $a \in S_{p(x_i)}$,
\[ \Delta((\lambda_i f_{n,x})^* f_{n_i, x_i})(x_i) = \overline{\lambda_i f_{n,x}(\frak s [\gamma_i])} f_{n_i, x_i}(\frak s[\gamma_i]) = | f_{n,x}(\frak s[\gamma_i])| f_{n_i, x_i}(\frak s[\gamma_i]) > 0.\]
Moreover, from the continuity of $\frak s$ and  $f_{n,x}$, and the fact that $f_{n,x}(\frak s[\gamma]) > 0$, we conclude that $\lim_i \lambda_i = 1$.  It follows that, for any open neighborhood $\mathcal O$ of $1 \in \T$, there exists $J$ so that whenever  $i \geq J$,  $\llbracket f_{n_i, x_i}, x_i \rrbracket = \llbracket  \lambda_ i f_{n,x}, x_i  \rrbracket \in Z(n, \mathcal O)$.  That is, $\psi_{\frak s}$ is continuous, as claimed.

Given the existence of a continuous section $\psi_{\frak s}: G_W \to \Sigma_W$, the formula 
\begin{equation}
\label{eq:cocycle-1}
C([n,\alpha_m(x)], [m, x]) = \psi_{\frak s}([n, \alpha_m(x)]) \psi_{\frak s}([m, x]) \psi_{\frak s}([nm, x])^{-1}\end{equation}
gives a 2-cocycle on $G_W$ which implements the twist $\Sigma_W$.  
To be precise, the right-hand side of \eqref{eq:cocycle-1} gives an element of the form $\llbracket b, \alpha_{nm}(x) =: y  \rrbracket \in \Sigma_W$, where $b$ is an element of the $\Gamma$-Cartan subalgebra $D = C^*_r(S, \omega)$.   (In fact, we may assume $b \in C_c(S, \omega)$ is supported in a bisection, as it is the product of elements which are so supported.)
Therefore (cf.~\cite[p. 47]{renault-cartan}), $\llbracket b, y \rrbracket$ can be identified with the element $(y(b)/|y(b)|, y) \in \T \times \widehat D$; and it is this element $y(b)/|y(b)| \in \T$ that we view as $C([n,\alpha_m(x)], [m, x])$, when we think of $C: G_W^{(2)} \to \T$. 

Recall that $\llbracket n, w \rrbracket^{-1} = \llbracket n^*, \alpha_n(w)\rrbracket$.  If the isomorphism $\Psi$ of Theorem \ref{Weyl groupoid} satisfies $\Psi([\gamma],\alpha_m(x)) = [n, \alpha_m(x)]$ and $\Psi([\eta], x) = [m, x]$, then 
\[  \alpha_{nm}(x) = \Psi(r(([\gamma], \alpha_m(x)), ([\eta], x))) =  \Psi([\gamma \eta] . x)  = [\gamma \eta].x.\]
  Consequently, 
we compute that 
\begin{align*}
\psi_{\frak s}& ([n, \alpha_m(x)])  \psi_{\frak s}([m, x]) \psi_{\frak s}([nm, x])^{-1}(\alpha_{nm}(x)) \\
& = \sum_{a \in S_{r(\gamma)}} \alpha_{nm}(x)(a) f_{n, \alpha_m(x)} f_{m, x} f_{nm, x}^*(a)\\
&= \sum_{a\in S_{r(\gamma)}} \left( x(\frak s([\gamma \eta])^{-1} a \frak s([\gamma \eta])) \omega(\frak s([\gamma \eta])^{-1}, a \frak s([\gamma \eta]) ) \overline{\omega(a \frak s([\gamma \eta]), \frak s([\gamma \eta])^{-1})} \right. \\ 
& \quad \times \left. \sum_{r(\xi) = r(\gamma)} (f_{n, \alpha_m(x)} f_{m, x} )( \xi) \overline{f_{nm, x}(a^{-1} \xi)} \overline{\omega( \xi^{-1}a,a^{-1} \xi)} \omega(\xi, \xi^{-1}a) \right) \\
 &= \sum_{a \in S_{r(\gamma)}} \left(x(\frak s([\gamma \eta])^{-1} a \frak s([\gamma \eta]))\omega(\frak s([\gamma \eta])^{-1}, a \frak s([\gamma \eta]) ) \overline{\omega(a \frak s([\gamma \eta]), \frak s([\gamma \eta])^{-1})} \right.\\
 & \quad \times  \left. \sum_{r(\xi) = r(\gamma)} \sum_{\xi = \kappa \lambda} f_{n, \alpha_m(x)}(\kappa) f_{m, x}(\lambda) \overline{f_{nm, x}(a^{-1} \xi)} \overline{\omega( \xi^{-1}a,a^{-1} \xi)} \right. \\
 & \quad \times \left. \omega(\xi, \xi^{-1}a) \omega(\kappa, \lambda) \right) .
\end{align*}
As everything in sight is supported on a bisection, there is at most one nonzero term in this sum; it occurs when $a^{-1} \xi = \frak s([\gamma \eta]),$ $\lambda = \frak s([\eta])$, and $\kappa = \xi \lambda^{-1} = a \frak s([\gamma \eta]) \frak s([\eta]) ^{-1} = \frak s([\gamma])$.  That is (writing, e.g., $f_{[\gamma]}:= f_{n, \alpha_n(x)}$),
\begin{align*}
\psi_{\frak s}& ([n, \alpha_m(x)])  \psi_{\frak s}([m, x]) \psi_{\frak s}([nm, x])^{-1}(\alpha_{nm}(x))  \\
& = x(\frak s([\gamma \eta])^{-1}\frak s([\gamma]) \frak s([\eta]) )\omega(\frak s([\gamma \eta])^{-1}, \frak s([\gamma]) \frak s([\eta]) ) \\
& \times \overline{\omega(\frak s([\gamma]) \frak s([\eta]), \frak s([\gamma \eta])^{-1})}  f_{[\gamma]}(\frak s[\gamma]) f_{[\eta]}(\frak s([\eta])   \overline{ f_{[\gamma \eta]}(\frak s[\gamma \eta])} \overline{\omega(\frak s([\gamma \eta])^{-1},\frak s([\gamma\eta])} \\
& \times\omega(\frak s([\gamma]) \frak s([ \eta]), \frak s([\gamma \eta])^{-1}) \omega(\frak s([\gamma]), \frak s([\eta]))\\
&= x(\frak s([\gamma \eta])^{-1}\frak s([\gamma]) \frak s([\eta]) ) \omega(\frak s([\gamma\eta])^{-1}, \frak s([\gamma]) \frak s([\eta])) \overline{\omega(\frak s([\gamma \eta])^{-1}, \frak s([\gamma \eta]))}  \\
& \times \omega(\frak s([\gamma], \frak s([\eta])) f_{[\gamma]}(\frak s[\gamma]) f_{[\eta]}(\frak s([\eta])   \overline{ f_{[\gamma \eta]}(\frak s[\gamma \eta])} .
\end{align*}
By construction, each of the $f$ terms is positive, and all the other terms lie in $\T$, so when we compute the associated circle element $C([n, \alpha_m(x)], [m, x]) = C(([\gamma], [\eta].x), ([\eta], x))$, we find that 
\begin{align*}
 C([n, &\alpha_m(x)], [m,x])  = C(([\gamma], [\eta].x), ([\eta], x))\\
 &  =  x(\frak s([\gamma \eta])^{-1}\frak s([\gamma]) \frak s([\eta]) ) \\
 & \quad \times \omega(\frak s([\gamma\eta])^{-1}, \frak s([\gamma]) \frak s([\eta])) \overline{\omega(\frak s([\gamma \eta])^{-1}, \frak s([\gamma \eta]))} \omega(\frak s([\gamma], \frak s([\eta])) \\
 &=  x(\frak s([\gamma \eta])^{-1}\frak s([\gamma]) \frak s([\eta]) ) \overline{\omega(\frak s([\gamma \eta]), \frak s([\gamma \eta])^{-1} \frak s([\gamma]) \frak s([ \eta]))} \omega(\frak s([\gamma], \frak s([\eta])).
 \end{align*}
 (The last equation follows from the cocycle identity and the fact (cf.~\cite[Lemma 2.1]{DGNRW}) that $\omega(\gamma, \gamma^{-1})  = \omega(\gamma^{-1}, \gamma)$ for all groupoid elements $\gamma$.)
\end{proof}

\section{Applications: Semidirect products}
\label{sec:applications}

In this section, we apply our results from Section \ref{sec:main-results} to identify conditions under which the twisted $C^*$-algebra of a semidirect product of groups can be ``untwisted'' into a groupoid $C^*$-algebra without a twist.  See Theorem \ref{thm:semidirect} and Corollary \ref{cor: semi} below.  In particular, the noncommutative solenoids of \cite{LP} satisfy our hypotheses. Indeed, that paper showcases the importance of having access to both the twisted and the untwisted descriptions of the $C^*$-algebras.  For example, Latr\'emoli\`ere and Packer use the untwisted transformation group picture in \cite[Theorem~3.4]{LP} to characterize simplicity of noncommutative solenoids, but use the twisted group picture in \cite[Theorem~3.7]{LP} to compute the $K$-theory for these algebras. 

Since Abelian groups, and extensions of them, are always amenable, if $H$ and $K$ are Abelian groups, their  semidirect product $H\rtimes K$ is also amenable.  Thus $C^*_r(H\rtimes K, \omega )\cong C^*(H\rtimes K, \omega)$ for any 2-cocycle $\omega$ on $H \rtimes K$,  so we drop  the subscript $r$ in what follows.

\begin{thm}
Suppose $H$ and $K$ are discrete Abelian groups and $\beta: K\to \text{Aut}(H)$ is an action of $K$ on $H$.   Let $G=H\rtimes_\beta K$ and   $\omega$ be a $2$-cocycle on $G$ such that $\omega|_H\equiv 1$.  Then $C_0(\widehat H) = C(\widehat H) \subseteq C^*(G, \omega)$ is $K$-Cartan and $C^*(G, \omega)\cong C^*(K\rtimes \widehat H, \omega|_K)$ where the action of $K$ on $\widehat H$ is given by 
\[
k.x(h)=\overline{\omega((0,-k),(0,k))}\omega((0,-k),(h,0))\omega((\beta_{-k}(h),-k),(0,k))x(\beta_{-k}(h)).
\]
\label{thm:semidirect}
\end{thm}

\begin{proof}
	
	We use Theorem \ref{gamma cartan exists}.
	Define $c: H\rtimes_\beta K\to K$ by $c(h,k)= k$.  Then $c$ is a group (hence groupoid) homomorphism with $c\inv (0)\cong H$.  As everything in sight has the discrete topology, $\Iso(c\inv(0))=H = \Iso(c^{-1}(0))^\circ$ is clopen in $H\rtimes_\beta K$.  Since $\omega|_{H} \equiv 1$, $H$ is maximal among Abelian subgroups of $\Iso(c^{-1}(0))^\circ$ on which $\omega$ is symmetric.  Moreover, $H = \Iso(c^{-1}(0))$ is immediately centralizing in $\Iso(c^{-1}(0))$, and  $H \leq H \rtimes_\beta K$ is normal by construction.  
%

	  Thus $C^*(H)\subseteq C^*(H\rtimes K)$ is $K$-Cartan.    By Theorem~\ref{Weyl groupoid}, the Weyl groupoid of $(C^*(H\rtimes K), C^*(H))$ is isomorphic to $((H\rtimes K)/H)\rtimes \widehat H\cong K\rtimes \widehat H$ with the action of $K$ on $\widehat H$ given by 
	
	\begin{align*}
		k.x(h)&=\overline{\omega((0,k),(0,-k))}\omega((0,-k),(h,0))\omega((0,-k)(h,0),(0,k))\\
		& \qquad \times x((0,-k)(h,0)(0,k))\\
		&=\overline{\omega((0,-k),(0,k))}\omega((0,-k),(h,0))\omega((\beta_{-k}(h),-k),(0,k))x(\beta_{-k}(h))
	\end{align*}
	as desired. 
	
		To see that $C^*(G, \omega) \cong C^*(K \rtimes \widehat H, \omega|_{K})$, we use Theorem \ref{thm:Weyl-twist}.  As all groups in question are discrete,  the section $\frak s: G/H \cong K \to G$ given by $\frak s(H, k) = (0,k)$ is continuous.  Therefore, Theorem \ref{thm:Weyl-twist} tells us that the Weyl twist over $K \rtimes \widehat H$ is given by the 2-cocycle 
	\begin{align*}
	C((k, k'.x), (k', x)) &= x((0, kk')^{-1} (0, k)( 0, k')) \overline{\omega((0, kk'), (0, kk')^{-1}(0, k)(0,k')) }\\
		& \qquad \times \omega ((0, k), (0, k'))\\
	&= x((0,0)) \overline{\omega((0, kk'), (0,0))} \omega((0, k), (0, k')) \\
	&= \omega((0, k),(0, k')). \qedhere
	\end{align*} 
\end{proof}

\begin{rmk} We can show that $C^*(H)$ is $K$-Cartan in $C^*(H\rtimes K, \omega)$ directly.  	
	First we show that the map $h\mapsto (h,0)$ induces an embedding of $C^*(H)$ in $C^*(H\rtimes K, \omega)$.    Recall that if $\Gamma$ is a discrete group,   $C^*(\Gamma ) = C^*(\{ \delta_\eta: \eta \in \Gamma\})$, where $\delta_\eta \delta_{\nu} = \delta_{\eta \nu}$ and $\delta_\eta^* = \delta_{\eta^{-1}} = \delta_\eta^{-1}$.  Thus, it suffices to show that $\delta_{(h,0)}\delta_{(h',0)}=\delta_{(h+h',0)}$.  This is a computation:
	\[
	\delta_{(h,0)}\delta_{(h',0)}=\delta_{(h,0)(h',0)}\omega((h,0),(h',0))=\delta_{(h +\beta_0(h'),0)}1=\delta_{(hh',0)}.
	\]
	As $H$ is discrete and Abelian, $C^*(H) $ is Cartan in itself.  Moreover, there is an evident grading of $C^*(H \rtimes K, \omega)$ by $K$, namely, $A_t := \text{span}\{ \delta_{(h,t)}: h \in H\},$ and the natural conditional expectation of $A$ onto $A_0$ is easily checked to be faithful.  Thus, to see that $C^*(H) \subseteq C^*(H\rtimes K, \omega)$ is $K$-Cartan, it suffices to see that, 
	for each $(h,k)\in C^*(H\rtimes K, \omega )$,  $\delta_{(h,k)}$ is a normalizer of $C^*(H)$ in $C^*(H\rtimes K)$.  Since  $\delta_{(h,k)}=\overline{\omega((h,0),(0,k))}\delta_{(h,0)}\delta_{(0,k)}$, we merely need check that $\delta_{(0,k)}^*\delta_{(h,0)}\delta_{(0,k)} \in C^*(H)$ for all $k \in K$.   This holds since
	\begin{align*}
		\delta_{(0,k)}^*\delta_{(h,0)}\delta_{(0,k)}& =\delta_{0,-k} \overline{\omega((0,k),(0,k)\inv)}\delta_{(h,0)(0,k)}\omega((h,0),(0,k))\\
		&= \overline{\omega((0,k),(0,-k))}\omega((h,0),(0,k))\delta_{(0,-k)}\delta_{(h,k)}\\
		&=\overline{\omega((0,k),(0,-k))}\omega((h,0),(0,k))\omega((0,-k),(h,k))\delta_{(0,-k)(h,k)}\\
		&=\overline{\omega((0,k),(0,-k))}\omega((h,0),(0,k))\omega((0,-k),(h,k))\delta_{(\beta_{-k}(h), -k k)}\\
		&=\overline{\omega((0,k),(0,-k))}\omega((h,0),(0,k))\omega((0,-k),(h,k))\delta_{(\beta_{-k}(h), 0)}\\
		& \in C^*(H).
	\end{align*}
\end{rmk}

\begin{cor} \label{cor: semi}
If in addition $\omega|_K\equiv 1$ then $C^*(G, \omega)\cong C^*(K\rtimes \widehat H)$ where the action of $K$ on $\widehat H$ is given by 
\[
k.x(h)=\omega((0,-k),(h,0))\omega((\beta_{-k}(h),-k),(0,k))x(\beta_{-k}(h)).
\]
\end{cor}

\begin{example}[Rotation algebras]
\label{ex:nc-tori}
	Fix $\theta \in \R$.
If $H=K=\Z$ and $\beta\equiv 1$  and $\omega_\theta ((h,k),(h',k'))=e^{2\pi i \theta k h'}$ then $C^*(H \rtimes K, \omega_\theta )$ is the rotation algebra $A_\theta$.  In this case, $\omega|_H(m,n)=\omega((m,0),(n,0))=e^{2\pi i \theta 0 n} \equiv 1$, and similarly $\omega|_K\equiv 1$, so   Corollary ~\ref{cor: semi} applies: 
\[ A_\theta \cong C^*(\Z^2, \omega_\theta) \cong C^*(\Z\rtimes \T), \]
where  the action of $\Z$ on $\T = \widehat \Z$ is given by   $(k.z)(h)=e^{2\pi i \theta (-hk)}z(h)$.  That is, 
\[ k . z = e^{-2\pi i \theta  k} z.\]
\end{example}

\begin{rmk}
In the setting of Corollary \ref{cor: semi}, one can detect simplicity of $C^*(G, \omega)$ from the action of $K$ on $\widehat H$:  \cite[Theorem 5.1]{BCFS} tells us that $C^*(G, \omega) \cong C^*(K  \rtimes \widehat H)$ will be simple iff $K \rtimes \widehat H$ is effective and minimal.
\end{rmk}


We now show how Theorem \ref{thm:semidirect} and Corollary \ref{cor: semi} apply  to the noncommutative solenoids of \cite{latremoliere-packer, LP}.
\begin{example}\label{ex: solenoid} For a natural number $N$, $\Z[\frac{1}{N}]$ can be obtained as the inductive limit $\Z\stackrel{m_N}{\to}\Z\stackrel{m_N}{\to}\Z\stackrel{m_N}{\to}\cdots$ where $m_N$ is the group homomorphism given by multiplication by $N$.   
	 Consequently, the dual of  $\Z[\frac{1}{N}]$ is given by
\[
\mathscr{S}_N:=\{(z_n)_{n\in \N}\in \T^\N: \forall n\in \N, z^N_{n+1}=z_n\};
\]
 the pairing is given by $\langle \frac{p}{N^i} ,(z_n)_n\rangle =z_i^p$  \cite[Proposition~1.2]{LP}.
 Observe that $\mathscr S_N$ is a compact group since $\Z[\frac{1}{N}]$ is discrete.

Fix a sequence $(\alpha_n)_n$ such that $\alpha_0\in [0,1)$ and for every $n$, there exists $k_n\in \{0,\ldots , N-1\}$ such that $N\alpha_{n+1}=\alpha_n+k_n$.  The  function  $\omega = \omega_\alpha: \Z[\frac{1}{N}] \times \Z[\frac{1}{N}] \to \T$ given by 
\begin{equation}
\label{eq:solenoid-cocycle}
 \omega\left(\left(\frac{p}{N^{h}}, \frac{q}{N^{k}}\right),\left(\frac{r}{N^{j}}, \frac{s}{N^{\ell}}\right)\right)=\exp\left(2\pi i \alpha_{h+\ell}ps\right)
 \end{equation}
 is a $2$-cocycle on $\Z[\frac{1}{N}]\times \Z[\frac{1}{N}]$ \cite[Theorem~2.1]{LP}, and by \cite[Theorem~2.3]{LP},  every $2$-cocycle on $\Z[\frac{1}{N}]$ is cohomologous to one of this form.  The resulting twisted group $C^*$-algebra $C^*( \Z[\frac{1}{N}]\times \Z[\frac{1}{N}], \omega)$ is called a noncommutative solenoid \cite[Definition~3.1]{LP} and is meant to generalize the  noncommutative tori.  
 
 As in the previous example, if $p=0$ or $s=0$ then $\omega=1$.  Thus, if we take $H=K=\Z[\frac{1}{N}]$ and $\beta\equiv 1$ in Corollary~\ref{cor: semi},  $\omega|_H\equiv 1\equiv \omega|_K$. This yields
$C^*( \Z[\frac{1}{N}]\times \Z[\frac{1}{N}], \omega)\cong C^*( \Z[\frac{1}{N}]\rtimes \mathscr{S}_N)\cong C\left(\mathscr{S}_N\right)\ltimes_\theta \Z[\frac{1}{N}]$ with the action of $\Z[\frac{1}{N}]$ on $\mathscr S_N$ given by 
$ {\frac{p}{N^\ell}.(z_n)_n (\frac{q}{N^j})=\left(\exp(2\pi i \alpha_{j+\ell}pq) z_j^q\right)}$.
That is, 
\[ \left( \frac{p}{N^\ell} . (z_n)\right)_n = \left( \exp(2\pi i \alpha_{n+\ell} p) z_n\right)_n.\]
This recovers
 \cite[Proposition~3.3]{LP}.   
\end{example}

\section{Reconstruction}

\label{sec:reconstruction}

Our motivating question in this section is: Can we recognize when a groupoid $H$ is of the form of our Weyl groupoids $G/S \rtimes \widehat S$? 

{Throughout this section, we will restrict our attention to the setting where the 2-cocycle $\omega$, when restricted to  $S$,  is trivial.  Indeed, some of the results in this section will place more restrictions on $\omega$.}

There are some properties that $H$ must  have: to begin with,  it must be \'etale,  and (assuming $\omega$ is trivial on $S$) its unit space must be a bundle $T$ of  Abelian groups, in such a way that the quotient $H/T$ is again a groupoid.  We show that these properties are nearly enough to guarantee that $H$ is of the form $G/S \rtimes \widehat S$ and to recover $G, S$ from $H$.

Suppose $H$ is a Hausdorff \'etale groupoid and $H^{(0)}=T$.   Suppose $p:T\to X$ is an open quotient map with $p\inv(x)$ a compact Abelian group for all $x$. Then we can consider $H$ as both a right and a left $T$-space, via the moment maps 
\[p_s:=p \circ s\quad \text{and}\quad p_r:=p \circ r.   \]   
We will use these actions to construct a new groupoid from $H$, $H/T\rectimes \widehat T$, in Theorem~\ref{prop:reconstruction-is-groupoid} below.  Then we show in Theorem~\ref{thm:uniqueness} that if $H=G/S\rtimes \widehat S$ is a Weyl groupoid with $\omega$ trivial on $G$, then there exists $\theta$ such that  $H/T\rectimes \widehat T\cong G$.   
 
 Our first goal is to define conditions that ensure that the quotient $H/T$ of $H$ under the action of $T$ is a groupoid.  This is not completely straightforward, because the group bundle (groupoid) structure of $T$ is not inherited from the groupoid structure of $H$.  Therefore, in order for $H/T$ to be a groupoid,  we need to assume that the right and left actions of $T$ on $H$ have some compatibility with the groupoid structure.  We are guided by the case that $H=G/S\rtimes \widehat S$ -- cf.~Proposition \ref{prop:weyl-satisfies} below.  As there are multiple actions associated to our constructions we will use different symbols for each, which we will introduce as the need arises.   To distinguish the left and right actions of $T$ on $H$ for $t\in T $ and $\eta\in H$ we denote
 \[
 t\ltact \eta \text{ (left action)}\quad \text{ and } \quad \eta\rtact t \text{ (right action)}.
 \]

To be precise, our assumptions on $H$ and $T = H\z$ are as follows.

\begin{assumption}
\label{assumption-on-actions}
{\color{white},}

\begin{enumerate}

\item Viewing $H$ as a $T$-space via $p_r$ and $p_s$, there exist commuting free and proper actions of $T$ on the left and right of $H$, such that
\begin{enumerate}

\item $t\ltact u=u \rtact t$ for all $u\in H^{(0)}, t\in T$.
\item for all $t\in T$, $\eta\in p_s\inv (p(t))$  and $\gamma\in p_r\inv(p(t))$ we have
 \[ r(t\ltact\gamma)=t\ltact r(\gamma)\quad\text{and}\quad  s(\eta\rtact t)=s(\eta)\rtact t\]
\end{enumerate}
\item There exist continuous maps $\lambda: H\mathbin{_{p_s}*_p}T\to T$ and $\rho: T\mathbin{_p*_{p_r}} H\to T$ such that for $\gamma,\eta\in H$ and $t\in T$, 
\begin{enumerate}
\item $\eta\rtact t=\lambda_\eta (t)\ltact \eta$
\item $t\ltact\eta=\eta\rtact\rho_\eta(t)$

\item $(\gamma\eta)\rtact t=(\gamma\rtact \lambda_\eta(t))(\eta\rtact t)$
\item $t\ltact(\gamma\eta)=(t\ltact\gamma)(\rho_\gamma(t)\ltact\eta)$
\item $(\eta\rtact t)\inv=
t\ltact \eta^{-1}$ 
\item $(t\ltact \eta)\inv=(\eta\inv)\rtact
t$
\end{enumerate}
\end{enumerate}
\end{assumption}
\begin{rmk}The following are immediate consequences of the conditions of Assumption \ref{assumption-on-actions}.
\label{rmk:consequences}
\begin{enumerate}
\item For all $\eta \in H$, $\lambda_\eta$ = $\rho_\eta^{-1}$ as we have
\[t\ltact\eta=\eta\rtact \rho_\eta(t)=\lambda_\eta(\rho_{\eta}(t))\ltact \eta,\]
and so from the freeness of the left action we conclude  that $t=\lambda_\eta(\rho_{\eta}(t))$ whenever $(t, \eta) \in T \mathbin{_p *_{p_r}} H $. Similarly $s=\rho_\eta(\lambda_\eta(s))$ whenever $(\eta, s) \in H\mathbin{_{p_s} *_p} T $. 
\item   Combining  point (1) above with points (a) and (b) of Assumption \ref{assumption-on-actions} (2), we conclude that $H/T = T\backslash H$.
\item A computation using the fact that the left and right actions commute implies that $\lambda_\eta(st) = \lambda_\eta(s) \lambda_\eta(t)$, and similarly for $\rho$.
\item Combining conditions 1(b) and 2(b) of Assumption \ref{assumption-on-actions}, we see that (writing $x_{\triv}$ for the identity operator on $T_x$)
\[ \rho_v = p(v)_{\triv} \quad \ \text{ for all }v \in H^{(0)},\]
and similarly for $\lambda$.
\end{enumerate}
\end{rmk}

As justification for Assumption \ref{assumption-on-actions}, the following proposition establishes that our motivating examples, the Weyl groupoids of $\Gamma$-Cartan pairs, satisfy Assumption \ref{assumption-on-actions}. 

\begin{prop}
\label{prop:weyl-satisfies}
Suppose $G$ is a Hausdorff \'etale groupoid and $S \leq \Iso(G)$ satisfies the hypotheses of Theorem~\ref{gamma cartan exists}, with a  2-cocycle $\omega$ which is trivial on $S$.
Set $H = G/S \rtimes \widehat S$, $T = \widehat S$, and let $p: \widehat S \to G\z$ denote the canonical projection map.  Consider the left action of $\widehat S$ on $H$ given by 
\[ \hat{t}\ltact([\gamma], \hat t') = ([\gamma], \Ad_{\gamma}(\hat t)\hat t')\]
 and the right action given by $([\gamma], \hat t')\rtact\hat t = ([\gamma], \hat t' \hat t)$.  For $([\gamma], \hat t') \in H$ and $\hat t \in T$, define  
 \[ \lambda_{([\gamma], \hat t')}(\hat t) =: \lambda_\gamma(\hat t) = \Ad_{\gamma\inv}(\hat t)\quad\text{and} \quad \rho_{([\gamma], \hat t')} (\hat t) =: \rho_\gamma(\hat t )  = \Ad_{\gamma}(\hat t).\]
Then $(H,T)$ with these maps satisfies the conditions of Assumption \ref{assumption-on-actions}.
\end{prop}
\begin{proof}
We will primarily discuss the left action and $\lambda$; the arguments for the right action and $\rho$ are similar but easier.
We first observe that our left action, as well as the maps $\rho, \lambda$, are well-defined.  This will follow once we establish that 
\[ \Ad_{\gamma}(\hat t) = a \mapsto \hat t(\gamma a \gamma^{-1})\]
 depends only on $[\gamma]$.  To that end, suppose  $[\gamma] = [\gamma']$, so that $ \gamma^{-1}  \gamma' \in S$.  As $S$ is normal, we also have $\gamma' \gamma^{-1} \in S$.  Thus,   if $\hat t \in \widehat S_{r(\gamma)}$, the fact that $\hat t: S_{r(\gamma)} \to \T$ is a homomorphism implies that, for all $a \in S_{s(\gamma)} = S_{s(\gamma')}$,
\[ \hat t(\gamma'  a  (\gamma')^{-1}) = \hat t (\gamma (\gamma ')^{-1}) \hat t(\gamma' a (\gamma')^{-1}) \hat t(\gamma' \gamma^{-1}) = \hat t (\gamma a \gamma^{-1}).\]

 The continuity of  multiplication in $G$ and in $\widehat S$ implies the continuity of the left and right actions as well as of $\lambda, \rho$.

To see that the left action is free, recall that the units in $\widehat S$ (viewed as a bundle of groups) are the trivial 1-dimensional representations of the fibre groups $S_u$.  Thus, suppose that $\hat t \ltact ([\gamma], \hat t') = ([\gamma], \hat t')$ for some $([\gamma], \hat t') \in G/S \rtimes \widehat S$.  That means that for all $a \in S_{p(\hat t')}$, $\hat t(\gamma a \gamma^{-1} ) \hat t'(a) = \hat t'(a)$.  As $\hat t'(a) \in \T$ for all $a$, we conclude that $\hat t(\gamma a \gamma^{-1}) = 1$ for all $a \in S_{p(\hat t')}$.  It follows that $\hat t$ is a unit.

To see that the left action is proper, suppose that $([\gamma_i], \hat t_i) \to ([\gamma], \hat t)$ and $\hat \tau_i {\ltact} ([\gamma_i], \hat t_i) \to ([\eta], \hat t')$ where $\hat\tau_i\in \widehat S$.  Then $[\eta] = [\gamma]$, and by the definition of the topology on $\widehat S$ \cite[Proposition 3.3]{DGN}, whenever $(a_i)_i \subseteq S$ satisfies $s(a_i) = p(\hat t_i)$ and $a_i \to a$, we have 
\[ \hat \tau_i(\gamma_i a_i \gamma_i^{-1}) \hat t_i(a_i) \to \hat t'(a) \quad \text{ and } \quad \hat t_i(a_i) \to \hat t(a).\]
We conclude that $\hat \tau_i(\gamma_i a_i \gamma_i^{-1}) \to \hat t'(a)\hat t(a) ^{-1}$, and we claim that it follows that the net $\{ \hat \tau_i\}_i$ is convergent. To that end, fix $\{b_i\}_i \subseteq S$ with $b_i \to b$ and $s(b_i) = p(\hat \tau_i) = r(\gamma_i).$  As we have assumed that $[\gamma_i]\to [\gamma]$, by passing to a subnet and relabeling, we obtain  a net $(c_i)_i \subseteq S$ with $c_i \gamma_i  \to \gamma$ (cf.~\cite[Proposition 2.12]{williams-gpoid} and \cite[Proposition 2.13.2]{fell-doran}).  Recall that $s(c_i) = r(c_i) =r(\gamma_i)$.

 Now, define $a_i := \gamma_i^{-1} c_i ^{-1} b_i c_i \gamma_i$.  Then $a_i \in S, s(a_i) = s(\gamma_i) $, and $a_i \to \gamma b \gamma^{-1}$. Moreover,  as each $\hat \tau_i$ is a homomorphism on $S$,
\[ \hat \tau_i(\gamma_i a_i \gamma_i^{-1} ) = \hat \tau_i( c_i) ^{-1}\hat \tau_i (b_i)  \hat \tau_i(c_i) = \hat \tau_i (b_i).\]
Since $\hat \tau_i(\gamma_i a_i \gamma_i^{-1}) \to \hat t'(a)\hat t(a) ^{-1}$, and $(b_i)_i$ was an arbitrary convergent net in $S$, 
 we conclude that (a subnet of) $\hat \tau_i$ converges.  That is, the left action is proper as claimed. 
 
As $(G/S \rtimes \widehat S)\z = \widehat S$, Properties (a) and (b) of (1) from Assumption~\ref{assumption-on-actions} follow quickly from the definitions.
 
Many of the assertions of Assumption \ref{assumption-on-actions} (2) are straightforward to check.  We have already seen that $\lambda$ and $\rho$ are well-defined and continuous.   Assertions (a) and (b) follow quickly from the definitions.  We check Assumption (d) (Assumption (c) is similar).

Suppose that $p(\hat t) = r(\gamma)$, so that $\hat t \ltact ([\gamma \eta], \hat t')$
 is  defined.  Note that $([\gamma \eta], \hat t') = ([\gamma], [\eta]. \hat t') ([\eta], \hat t')$, where $[\eta].\hat t'$ denotes the action of $G/S$ on $\widehat S$ from Equation \eqref{eq:Weyl-gpoid-action}.  We begin by computing the right-hand side of 2(d).
\begin{align*}
\hat t \ltact ([\gamma], [\eta].(\hat t'))&=([\gamma], \Ad_\gamma(\hat t) [\eta].\hat t') \\
 \rho_\gamma(\hat t) \ltact ([\eta], \hat t')  &=  \Ad_\gamma( \hat t)\ltact ([\eta], \hat t')=([\eta], \Ad_{\gamma\eta}(\hat t)\hat t').
\end{align*}
One checks that, for any $a \in S_{r(\eta)}$, 
\begin{align*} 
r([\eta], \Ad_{\gamma \eta}(\hat t) \hat t') (a) & = [\eta].(\Ad_{\gamma\eta}(\hat t) \hat t'))(a) = \omega(\eta^{-1}, a \eta) \overline{\omega(a\eta, \eta^{-1})} \hat t(\gamma a \gamma^{-1}) \hat t'(\eta^{-1} a \eta) \\
&= \hat t(\gamma a \gamma^{-1} )\alpha_\eta(\hat t')(a)=\Ad_\gamma(\hat t) [\eta].\hat t'(a) = s( [\gamma], \Ad_\gamma(\hat t)[\eta].\hat t'),
\end{align*} 
so $( \hat t \ltact ([\gamma], [\eta].(\hat t')) ,  \rho_\gamma(\hat t) \ltact ([\eta], \hat t') ) \in H\2$. Indeed, their product is
\begin{align*}
\hat t \ltact ([\gamma], [\eta].(\hat t')) \rho_\gamma(\hat t) \ltact ([\eta], \hat t')&=([\gamma], \Ad_\gamma(\hat t) \alpha_\eta(\hat t'))([\eta], \Ad_{\gamma\eta}(\hat t)\hat t')\\
&=([\gamma][\eta],\Ad_{\gamma\eta}(\hat t)\hat t')=\hat t\ltact ([\gamma][\eta],\hat t')
\end{align*}
 as desired. 
 
 To see (f), we compute:
 \begin{align*}
 (\hat t\ltact ([\eta],\hat t'))\inv&=([\eta],\Ad_\eta(\hat t) \hat t')\inv = ([\eta^{-1}], [\eta].(\Ad_\eta(\hat t) \hat t')) \\ 
&= ([\eta^{-1}], \hat t ( [\eta]. \hat t'))
 =([\eta\inv], [\eta].\hat t')\rtact \hat  t=(([\eta],\hat t')\inv) \rtact \hat t.\end{align*}   The computation for (e) is similar.
\end{proof}

Indeed, Assumption \ref{assumption-on-actions} encapsulates nearly all the properties of the Weyl groupoid $G/S \rtimes \widehat S$: Theorems \ref{thm:reconstruction-1} and \ref{thm:uniqueness} explain  how to recover the initial data $(G, S)$ from the Weyl groupoid $H = G/S \rtimes \widehat S$.  

Our first goal is to construct a groupoid $H/T\rectimes \widehat T$ from the data in Assumption~\ref{assumption-on-actions}; this groupoid will recover the groupoid $G$ if we start with $H = G/S \rtimes \widehat S$.  We do this in a series of lemmas. 

Let $[\eta]$ be the equivalence class of $\eta$ in $H/T = T \backslash H$.  Our first step is to show that $H/T$ inherits an \'etale groupoid structure from $H$.  The next two lemmas will be useful in establishing the topology and algebraic structure on $H/T$.

\begin{lemma}\label{lemma: quotient bisection} Suppose that $\gamma\in H$ and $B$ is a bisection containing $\gamma$.  Then there exists an open neighborhood $W\subseteq B$ of $\gamma$ such that the range and source maps are injective on $[W]$.\end{lemma}

\begin{proof} 
We will first establish the existence of a neighborhood $W_r$ of $\gamma$ so that the range map is injective on $[W_r]$.  Analogously, we find a neighborhood $W_s$ such that the source map is injective on $[W_s]$. We then take $W := W_r \cap W_s$. 

To establish the existence of $W_r$, for a contradiction, suppose that it does not exist.  Then for every open neighborhood $W \subseteq B$ of   $\gamma$, we have $\eta^W, \gamma^W$ in $W$ such that $[\eta^W]\neq [\gamma^W]$ and $[r(\eta^W)]=[r(\gamma^W)]$.  Since $[r(\eta^W)]=[r(\gamma^W)]$ there exists $t_W\in T$ such that $t_W\ltact r(\eta^W)=r(\gamma^W)$.  {Then, $(t_W, \eta^W) \in T \mathbin{_p *_{p_r}} H$, and so  $t_W \ltact \eta^W$ makes sense.  We have $r(t_W \ltact \eta^W) = t_W \ltact r(\eta^W) = r(\gamma^W)$, so the fact that $B \supseteq W$ is a bisection implies that if $t_W \ltact \eta^W \in W$ then $t_W \ltact \eta^W = \gamma^W.$ 
}
 By definition, $\eta^W\to \gamma$ so $r(\eta^W)\to r(\gamma)$ and also $ r(\gamma^W) = t_W\ltact r(\eta^W)\to r(\gamma)$.  Since the left action is proper, we have a convergent subnet $t_{W_i}$ of the $t_W$'s. Since the action is free, this subnet must converge to the unit $ p_r(\gamma) \in X$.  
 Thus 
 \[ t_{W_i}\ltact\eta^{W_i}\to  p_r(\gamma) \ltact \gamma = \gamma.\]    It then follows that $t_{W_i}\ltact \eta^{W_i} \in B$ eventually.  
As $B$ is a bisection, $\gamma^{W_i} \in B$ for all $i$, and $r(t_{W_i} \ltact \eta^{W_i}) = r(\gamma^{W_i})$, we conclude that  $t_{W_i} \ltact \eta^{W_i} = \gamma^{W_i}$ for large enough $i$.   This contradicts our assumption that  $[\eta^W]\neq [\gamma^W]$ for all $W$.    Therefore, there must exist a neighborhood $W_r \subseteq B$ of $\gamma$ such that $r_{H/T}|_{[W_r]}$ is injective.  Repeat this process replacing $r$ with $s$ to get an open neighborhood $W_s$ of $\gamma$ such that $s_{H/T}|_{[W_s]}$ is injective.  Now take $W=W_r\cap W_s$.
 \end{proof}

\begin{lemma} \label{lem equiv are bisections} Suppose  $[\gamma], [\eta]\in H/T$.  If there exists $\gamma'\in [\gamma]$ and $\eta'\in [\eta]$ such that $(\gamma',\eta')\in H^{(2)}$, then for every $\gamma''\in [\gamma]$ there exists a unique $\eta''\in [\eta]$ such that $(\gamma'',\eta'')\in H^{(2)}$.   
\end{lemma}

\begin{proof}
Suppose $\gamma''\in [\gamma]$.  Then there exists $t\in T$ such that $\gamma''=\gamma'\rtact t $.  By hypothesis, $s(\gamma') = r(\eta')$, so $(t, \eta') \in T \mathbin{_p *_{p_r}} H$ and we can  consider $\eta'' = t\ltact \eta' =\eta'\rtact\rho_{\eta'}(t)$.    By Assumption \ref{assumption-on-actions}(1),
\[ s(\gamma'')=s(\gamma')\rtact t=t\ltact r(\eta') =r(t\ltact \eta') 
=r(\eta'').\]
   Thus $(\gamma'',\eta'')\in H^{(2)}$.    Now if $(\gamma'', \eta''')\in H^{(2)}$ with $\eta'''\in [\eta]$ then there is $\tau\in T$ such that $\eta'''=\eta'\rtact \tau$ and  
   \[ t\ltact r(\eta')=r(\eta'')=s(\gamma'')=r(\eta''')=r(\lambda_{\eta'}(\tau)\ltact \eta')=\lambda_{\eta'}(\tau)\ltact r(\eta').\]
     The freeness of the left action implies that  $t=\lambda_{\eta'}(\tau)$ and so $\tau=\rho_{\eta'}(t)$.  In other words, $\eta''' = \eta' \rtact \rho_{\eta'}(t) = t\ltact \eta' =\eta''$ as desired.
\end{proof}

Now we will define a groupoid structure on $H/T$.

 \begin{prop}\label{prop: quotient groupoid} Suppose $(H, T)$ satisfies Assumption~\ref{assumption-on-actions}.  Then $H/T$ is a Hausdorff \'etale groupoid under the quotient topology  with \[(H/T)^{(2)}=\{([\gamma], [\eta]): \exists \gamma'\in [\gamma], \eta'\in [\eta] \text{ such that } (\gamma',\eta')\in H^{(2)}\}.\]
Composition and inversion are given by \begin{equation}\label{eq:composition} [\gamma][\eta] =  \{ \xi :   \xi = \gamma' \eta' \text{ for some } \gamma' \in [\gamma], \eta'\in [\eta]\}, \qquad  [\gamma]\inv =[\gamma\inv].\end{equation}
Moreover, $r_{H/T}([\gamma])=[r_H(\gamma)]$ and $s_{H/T}([\gamma])=[s_H(\gamma)]$.
\end{prop}

\begin{proof} 
We check that the composition law is well-defined. Let $ ([\gamma], [\eta])\in (H/T)^{(2)}$. We need to show that $[\gamma][\eta]$ is an equivalence class and further that $[\gamma][\gamma]\inv [\gamma]=[\gamma]$.  To see that $[\gamma][\eta]$ is an equivalence class, pick $\gamma',\eta'$ as in the definition of $(H/T)^{(2)}$. We claim $[\gamma][\eta]=[\gamma'\eta']$. To see this note that if $\zeta\in [\gamma'\eta']$ then there exists a $t$ such that $\zeta=(\gamma'\eta')\rtact t=(\gamma'\rtact \lambda_{\eta'}(t))(\eta'\rtact t)\in [\gamma][\eta]$.  That is,  $[\gamma'\eta']\subseteq [\gamma][\eta]$.  Now suppose $\zeta\in [\gamma][\eta]$.  Then there exist $\gamma''\in [\gamma]$ and $\eta''\in [\eta]$ such that $\gamma''\eta''=\zeta$.  As $[\eta] = [\eta']$, there exists $t\in T$ such that $\eta'\rtact t=\eta''$.  Thus $r(\eta'')=\lambda_{\eta'}(t)\ltact r(\eta')=s(\gamma')\rtact \lambda_{\eta'}(t)$. This implies that $(\gamma'\rtact\lambda_{\eta'}(t),\eta'')\in H^{(2)}$.  Consequently, the uniqueness statement of Lemma~\ref{lem equiv are bisections} yields $\gamma''=\gamma'\rtact\lambda_{\eta'}(t)$. In other words, for any $\zeta \in [\gamma] [\eta]$,
\[ \zeta=\gamma''\eta''=(\gamma'\rtact\lambda_{\eta'}(t))(\eta'\rtact t)=(\gamma'\eta')\rtact t \in [\gamma' \eta'].\]
 Thus $[\gamma][\eta]$ is an equivalence class as desired. 

Next we need to check that $[\gamma][\gamma]\inv [\gamma]=[\gamma]$.  First recall from Assumption~\ref{assumption-on-actions}(2)(e) that  $(\eta\rtact t)^{-1} =  t \ltact \eta^{-1} = \eta^{-1}\rtact \rho_{\eta^{-1}}(t)$. So  $[\gamma]^{-1} = [\gamma\inv]$ is well-defined and $[\gamma][\gamma]\inv [\gamma]=[\gamma]$ follows.  
We have  $r_{H/T}([\gamma])=[\gamma][\gamma\inv]=[\gamma\gamma\inv]=[r_H(\gamma)] $ and $s_{H/T}([\gamma]) = [\gamma]^{-1}[\gamma] = [s_H(\gamma)].$

Multiplication and inversion are continuous in $H/T$ since these operations are inherited from $H$ and the quotient map is continuous. 
Since the action of $T$ is proper, \cite[Proposition 2.18]{williams-gpoid} implies that  $H/T$ is Hausdorff.  
Similarly, 
$r_{H/T}$ is continuous, being the composition of continuous maps.  

To show that $H/T$ is \'etale we need to show that every element is contained in an open bisection.  Let $[\gamma]\in H/T$.   Lemma~\ref{lemma: quotient bisection} ensures the existence of a bisection $W$ in $H$ containing $\gamma$ such that the range and source maps are injective when restricted to $[W]$; the continuity of $r|_{[W]}, s|_{[W]}$ follows from the previous paragraph.  It remains to show these maps are 
open.    Let $U\subseteq [W]$ be open in $[W]$. Since $W$ is open in $H$ and the quotient map is open, $[W]$ is open in $H/T$, and hence $U$ is open in $H/T$.    Observe that  $r_{H/T}(U)=q(r_H(q\inv(U))$ is open, since $q$ is continuous and open and  $r_H$ is open.   That is, $r_{H/T}|_{[W]}$ is a homeomorphism as desired.  A symmetric argument will establish that $s_{H/T}|_{[W]}$ is a homeomorphism as well.
\end{proof}

Our next goal in the construction of $H/T\rectimes \widehat T$ is to establish an action of $H/T$ on $\widehat T$; we do this in Lemma~\ref{lem:T-hat-action} below.  First, however,  we need to establish that $\gamma\mapsto \rho_\gamma$ is a homomorphism.  

\begin{lemma}
For any $(\gamma, \eta) \in H^{(2)}$, 
$\lambda_{\gamma\eta}=\lambda_\gamma\circ\lambda_\eta$ and similarly $\rho_{\gamma\eta}=\rho_\eta\circ\rho_\gamma$.
\label{lem:lambda-multiplicative}
\end{lemma}

\begin{proof}
Consider \[\lambda_{\gamma\eta}(t)\ltact\gamma\eta=(\gamma\eta)\rtact t=(\gamma\rtact \lambda_\eta(t))(\eta\rtact t)=(\lambda_\gamma(\lambda_\eta(t))\ltact \gamma)(\lambda_\eta(t)\ltact \eta).\]
Since $\rho_\gamma(\lambda_\gamma(\lambda_\eta(t)))=\lambda_\eta(t)$, Assumption 5.1(2)(d) implies that the above is equal to  $(\lambda_\gamma(\lambda_\eta(t)))\ltact (\gamma\eta)$.  By the freeness of the left action, we conclude that  $\lambda_{\gamma\eta}(t)=\lambda_\gamma(\lambda_\eta(t))$ as desired.  As $\rho_\xi=(\lambda_\xi)\inv$, taking the inverse of the first equation gives  the second.\end{proof}

The following Lemma establishes that we have a well-defined action of $H/T$ on $\widehat T$.  

\begin{lemma} 
We have  an action of $H$ on $\widehat T$, given for $\chi \in \widehat T_{p_s(\gamma)}$ by $\gamma\hact\chi=\chi\circ\rho_\gamma$. Furthermore, for any $t \in T_{p_s(\gamma)}$ we have  $(\gamma\rtact t)\hact\chi=\gamma\hact\chi$.
Consequently, this action  descends to a continuous action of $H/T$ on $\widehat T$.  Moreover, $[\gamma]\hact(\chi_1 \chi_2) =( [\gamma]\hact\chi_1)([\gamma]\hact\chi_2)$.
\label{lem:T-hat-action}
\end{lemma}

\begin{proof}
 The fact that $\gamma \hact \chi = \chi \circ \rho_\gamma$ is an action of $H$ on $T$ follows from the fact that $\rho$ is a homomorphism, and that (Remark \ref{rmk:consequences}(4))  $\rho_v = p(v)_{\triv}$ is the identity operator.  We also have $p(\gamma \hact \chi) = p_r(\gamma)$,  since $\gamma \hact \chi(t) = \chi \circ \rho_\gamma(t)$ makes sense precisely when $t \in T_{p_r(\gamma)}$.

To see that $\hact$ descends to an action of $H/T$, we will  show that $\rho_{\gamma \rtact t} = \rho_{\gamma}$.  Thus, fix  $\gamma\in H$ and $t \in T_{p_s(\gamma)},\tau\in T_{p_r(\gamma \rtact t)}.$ 
As $T$ is a bundle of Abelian groups and $\rho$ is multiplicative, 
\begin{align*} \gamma\rtact ( \rho_{\gamma\rtact t}(\tau) t) &= (\gamma\rtact t)\rtact \rho_{\gamma\rtact t}(\tau)=\tau\ltact(\gamma\rtact t)=\tau\ltact(\lambda_{\gamma}(t)\ltact\gamma)=(\tau\lambda_{\gamma}(t))\ltact\gamma\\
&=\gamma\rtact(\rho_{\gamma}(\tau)t).\end{align*}
 Thus  the  freeness of the right action yields $\rho_{\gamma\rtact t} (\tau)t=\rho_{\gamma}(\tau)t$.  In particular, this yields $p_r(\gamma \rtact t)  = p_r(\gamma)$ for all $t \in T_{p_s(\gamma)}$.
 
  The fact that $t$ was arbitrary now implies that $\rho_{\gamma\rtact t}(\tau)=\rho_{\gamma}(\tau)$ for all $t\in T_{p_s(\gamma)}$ and all $\tau \in T_{p_r(\gamma)}$.  In other words, $\rho_{\gamma \rtact t} = \rho_\gamma$.  It  follows that $\hact$ descends to an action of $H/T$ on $\widehat T$, as claimed.

The fact that the given formula yields a continuous groupoid action of $H$, and hence of $H/T$, follows since $\rho$ is multiplicative
(Lemma \ref{lem:lambda-multiplicative}) and continuous (by assumption).  The last statement is immediate from the definitions, and the fact that multiplication in $\widehat T$ is pointwise within fibres.
\end{proof}
%

At this point we can construct a new groupoid from $(H, T)$.

\begin{prop} Suppose $H$ satisfies Assumption~\ref{assumption-on-actions}: in particular $H\z=T$ and $p: T\to X $ is an Abelian group bundle.  Suppose further that 
  $\theta: (H/T)*(H/T)\to \widehat T$ is continuous with $\theta([\gamma],[\eta])\in \widehat T([s(\eta)])$ and satisfies an analogue of the cocycle condition: for any $\gamma \in H$ we have  $\theta([r(\gamma)], [\gamma]) = [s(\gamma)]_{\triv}$, and whenever $([\gamma], [\eta]), ([\eta], [\zeta]) \in (H/T)\2$, 
\begin{equation}
\label{eq:omega-for-recovering-G}
([\zeta]\inv\hact\theta([\gamma],[\eta]))\theta([\gamma][\eta],[\zeta])=\theta([\gamma],[\eta][\zeta])\theta([\eta],[\zeta]).
\end{equation} 
Then 
 \[H/T\rectimes \widehat T:=\{([\gamma],\chi)\in H/T\times \widehat T: s([\gamma])=p(\chi)\}\]
 has a locally compact \'etale groupoid structure under the restriction of the product topology where 
 \begin{align*}
 (H/T\rectimes \widehat T)^{(2)}&=\{(([\gamma],\chi),([\eta],\nu)): s([\gamma])=r([\eta])\}\\
\intertext{ with composition and inversion given by}
([\gamma],\chi)([\eta],\nu)&=([\gamma][\eta], \theta([\gamma],[\eta])([\eta]\inv\hact\chi)\nu), \\
([\gamma],\chi)\inv&=([\gamma]\inv, \theta([\gamma],[\gamma]\inv)\inv([\gamma]\hact\chi\inv)),\\
\intertext{and range and source maps}
r([\gamma], \chi) &= ([r(\gamma)], [r(\gamma)]_{\triv})  \quad \text{ and } s([\gamma], \chi) = ([s(\gamma)], [s(\gamma)]_{\triv}).
\end{align*}
In particular, 
 $ (H/T\rectimes \widehat T)\z$ can be identified with $X$.

\label{prop:reconstruction-is-groupoid}
\end{prop}

\begin{rmk}
Observe that $H/T \rectimes \widehat T$ is not a (twisted) action groupoid.  First, $(H/T \rectimes \widehat T)^{(0)}\not= \widehat T.$  Also, $\theta$ takes values in $\widehat T$, not in $\T$.
\end{rmk}

\begin{rmk}
One can take $\theta$ to be trivial in Propostion~\ref{prop:reconstruction-is-groupoid} and still obtain a groupoid.  However in our reconstruction theorem below (Theorem~\ref{thm:uniqueness}) we need to allow for nontrivial $\theta$.  
\end{rmk}

\begin{proof}[Proof of Proposition~\ref{prop:reconstruction-is-groupoid}:]

The topology on $H/T \rectimes \widehat T$ is inherited from that of $H/T$ and $\widehat T$ via the (fibre-)product topology. The continuity of the action of $H/T$ on $\widehat T$, together with the continuity of $\theta$ and the continuity of multiplication and inversion in $H/T$ and $\widehat T$, implies that the groupoid operations on $H/T \rectimes \widehat T$ are continuous, being inherited from continuous maps.

To show $H/T \rectimes \widehat T$ is a groupoid we need to show that multiplication is associative and $([\gamma],\chi)([\gamma],\chi)\inv([\gamma],\chi)=([\gamma],\chi)$ for all $([\gamma], \chi) \in H/T\rectimes \widehat T$.

For the first assertion, since $\widehat T$ is a bundle of Abelian groups and $\hact$ is an action, we compute
\begin{align*}\left[([\gamma],\chi)([\eta],\nu)\right]&([\zeta], \mu)=\left([\gamma][\eta], \theta([\gamma],[\eta])([\eta]\inv\hact\chi)\nu\right)([\zeta], \mu)\\
&=\left(([\gamma][\eta])[\zeta],\theta([\gamma][\eta],[\zeta])([\zeta]\inv\hact( \theta([\gamma],[\eta])([\eta]\inv\hact\chi)\nu))\mu\right)\\
&=\left([\gamma]([\eta][\zeta]),[\zeta]\inv\hact( \theta([\gamma],[\eta]))\theta([\gamma][\eta],[\zeta])(([\eta ][\zeta])\inv\hact\chi)([\zeta]\inv\hact\nu)\mu)\right)\\
&=\left([\gamma]([\eta][\zeta]),\theta([\gamma],[\eta][\zeta])\theta([\eta],[\zeta])(([\eta][\zeta])\inv\hact\chi)([\zeta]\inv\hact\nu)\mu)\right)\\
&=\left([\gamma],\chi\right)\left([\eta][\zeta], \theta([\eta],[\zeta])([\zeta]\inv\hact\nu)\mu)\right)\\
&=([\gamma],\chi)\left[([\eta],\nu)([\zeta], \mu)\right].
\end{align*}

For the second assertion, if $p_r(\gamma) = x$, we have
\begin{align*}
([\gamma],\chi)&([\gamma],\chi)\inv([\gamma],\chi) =([\gamma], \chi)([\gamma]\inv,  \theta([\gamma],[\gamma]\inv)\inv( [\gamma]\hact\chi\inv))([\gamma], \chi)\\
&=([\gamma][\gamma]\inv,\theta([\gamma],[\gamma]\inv)\theta([\gamma],[\gamma]\inv)\inv ([\gamma]\hact\chi)( [\gamma]\hact\chi\inv))([\gamma], \chi)\\
&=  ([r(\gamma)], [\gamma]\hact [s(\gamma)]_{\triv})( [\gamma], \chi) = 
([r(\gamma)],x_{\triv})([\gamma], \chi)) = ([\gamma], \chi)
\end{align*}
as desired. \end{proof}

The following two Theorems show that, first, if we start with a pair $(H, T)$ arising as in Proposition \ref{prop:weyl-satisfies}, so that $H = G/S \rtimes \widehat S$ is the Weyl groupoid of a dynamical $\Gamma$-Cartan pair, then the construction of Proposition \ref{prop:reconstruction-is-groupoid} recovers $G$ and $S$ from the Weyl groupoid $H = G/S \rtimes S$.  Then,  Theorem \ref{thm:reconstruction-1} identifies  additional hypotheses  on $(H, T)$ which  guarantee that we obtain a pair $(G, S)$ which satisfies the hypotheses of Theorem \ref{gamma cartan exists}. 

We remark that our hypotheses in Theorem \ref{thm:reconstruction-1} are sufficient but almost certainly not necessary.  For one thing, our construction does not yield a 2-cocycle on $G$.  Moreover, in Theorem \ref{gamma cartan exists}, $S$ need not equal $\Iso(c^{-1}(0))^\circ$.  However, the hypotheses of Theorem \ref{thm:reconstruction-1} will always yield $S = \Iso(c^{-1}(0))^\circ$.

\begin{thm} 
Let $(G, S, c)$ satisfy the hypotheses of Theorem \ref{gamma cartan exists}, with a trivial 2-cocycle $\omega$. Let $\frak{s}:G/S\to G$ be a continuous section with $\frak s([u]) = u$ for all $u \in G\z$, and $H=G/S\rtimes \widehat S$.   Let $T=\widehat S$ and $\theta([\gamma],[\eta])=\frak{s}([\gamma\eta])\inv\frak{s}([\gamma])\frak{s}([\eta])$.   Then the map
$\Phi:H/T\rectimes \widehat T \to G$ given by 
\[ \Phi([\gamma],t) =  \frak{s}([\gamma])t\]
 is an isomorphism of topological groupoids.
 
 Moreover, $c$ induces a cocycle $\overline c : H \to \Gamma$ by $\overline c([\gamma], \hat t) = c(\gamma),$
 and {$c \circ \Phi = \overline c$.}
\label{thm:uniqueness}
\end{thm}

\begin{proof}
We first note that since $H=G/S\rtimes \widehat S$ and $T=\widehat S$ we have $H/T\cong G/S$ and $\widehat T\cong S$.  Thus, if $([\gamma], t ) \in H/T \rectimes \widehat T$, then $p(t) = s([\gamma]) = [s(\gamma)]$, so  $\frak{s}([\gamma])t$ makes sense as an element of $G$.  Indeed, the action of $S$ on $G$ preserves units.  Thus, if $([\gamma], [\eta]) \in (H/T)\2 = (G/S)\2,$ we have  $s(\gamma) = r(\eta)$ and consequently $[\gamma] [\eta] = [\gamma \eta]$.

We now check that the given formula for $\theta$ satisfies Equation \eqref{eq:omega-for-recovering-G}.  Using the definition for $\rho_\gamma$ from Proposition \ref{prop:weyl-satisfies},  we see that the action of $(H/T = G/S)$ on $\widehat T = S$ is given by 
\[ ([\gamma] \hact s) (t) = t(\tilde \gamma s (\tilde \gamma)^{-1}),\]
for any $\tilde \gamma \in [\gamma]$.
In other words, written as an action of $G/S$ on $S$, 
\[ [\gamma] \hact s = \tilde \gamma s \tilde \gamma^{-1}.\]
In the computation that follows, we will use the representative 
\[ \tilde \zeta := \frak s([\gamma \eta])^{-1} \frak s ([\gamma \eta \zeta]) \]
of $[\zeta]$ to compute the  action of $[\zeta]^{-1}$ on $\theta([\gamma], [\eta])$.
Together with the fact that $[\gamma] [\eta] = [\gamma \eta]$, we obtain that 
\begin{align*}
[\zeta]^{-1} \hact \theta([\gamma], [\eta]) \theta([\gamma][\eta], [\zeta]) & =\frak s([\gamma \eta \zeta])^{-1}\frak s([\gamma \eta])  \frak s([\gamma \eta])^{-1} \frak s([\gamma]) \frak s ([\eta])  \frak s([\gamma \eta])^{-1} \frak s([\gamma \eta \zeta]) \\
& \quad \times  \frak s([\gamma \eta \zeta])^{-1} \frak s([\gamma \eta]) \frak s([\zeta]) \\ 
&= \frak s([\gamma \eta \zeta])^{-1} \frak s ([\gamma]) \frak s([\eta]) \frak s([\zeta]), \\
\text{while } \ \theta([\gamma], [\eta] [\zeta]) \theta([\eta], [\zeta]) &= \frak s([\gamma \eta \zeta])^{-1} \frak s([\gamma]) \frak s([\eta \zeta]) \frak s([\eta \zeta])^{-1} \frak s([\eta]) \frak s([\zeta]) \\
&= \frak s([\gamma\eta \zeta])^{-1} \frak s([\gamma]) \frak s([\eta]) \frak s([\zeta]).
\end{align*}

Furthermore, since $\frak s([u]) = u$ for any $u \in G\z$,
\[ \theta([r(\gamma)], [\gamma]) = \frak s([\gamma])^{-1} \frak s([r(\gamma)]) \frak s([\gamma]) = s(\frak s [(\gamma)]) = s(\gamma).\]
Viewed as an element of $\widehat S$, the unit $s(\gamma)$ is the identity element.  That is, $\theta([r(\gamma)], [\gamma]) = s(\gamma)_{\triv}$ as desired.


We now check that $\Phi$ is a homomorphism: 
\begin{align*} \Phi(([\gamma],t)([\eta],t'))&=\Phi([\gamma\eta],\theta([\gamma],[\eta])([\eta]\inv \hact t)t')\\
&=\frak{s}([\gamma\eta])\frak{s}([\gamma\eta])\inv \frak{s}([\gamma])\frak{s}([\eta])([\eta]\inv \hact t)t'\\
&=\frak{s}([\gamma])\frak{s}([\eta])([\eta]\inv\hact t)t'\\
&=\frak{s}([\gamma])\frak{s}([\eta])\frak{s}([\eta])\inv t\frak{s}([\eta])t'\\
&=\frak{s}([\gamma])t\frak{s}([\eta])t'\\
&=\Phi([\gamma],t)\Phi([\eta],t').
\end{align*}

Now $\Phi$ is continuous since multiplication and $\frak{s}$ are continuous.  It thus suffices to see that $\Phi$ has a continuous inverse.  For $\gamma\in G$, define $\sigma(\gamma)=\frak{s}([\gamma])\inv \gamma$.  Since $\frak s$ is continuous and multiplication and inversion are continuous, $\sigma$ is continuous and takes values in $S=\widehat T$.  The map $\gamma\mapsto ([\gamma], \sigma(\gamma))$ is the inverse for $\Phi$. Indeed, 
\[
\Phi([\gamma],\sigma(\gamma))=\frak{s}([\gamma])\frak{s}([\gamma])\inv \gamma=\gamma. 
\]

For the final assertion, we first observe that $\bar c$ is well defined, since $c$ is a homomorphism and $S \subseteq c^{-1}(0)$.  In particular, for any $\gamma \in G$ and $t \in \widehat T = S$, 
\[ c(\gamma) = c(\frak s ([\gamma])) = c(\frak s([\gamma]) t), \]
and so $\bar c = c \circ \Phi$ as claimed.  As $c$ and $\Phi$ are multiplicative, it follows that $\bar c$ is indeed a cocycle.
\end{proof}

For our final result, we need to rephrase the notion of ``immediately centralizing'' (Definition \ref{def:imm-cent}) in terms of actions rather than subgroupoids.   

\begin{defn}
\label{def:imm-cent-action}
Let  $K$ be a group bundle and $\Gamma$ be a groupoid with a continuous open surjection $p: K \to \Gamma\z$.
A (left) action of  $\Gamma$ on $K$ is an {\em immediately centralizing action} if, for all integers $k \geq 1$ and all $g \in \Iso(\Gamma)$, 
\begin{equation}
\left[ (\forall \chi \in K \ni p(\chi) = s(g))(\exists 1 \leq n \leq k)( (g.\chi )^n = \chi^n ) \right] \implies g.\chi = \chi \ \forall \chi .
\label{eq:imm-cent-action}
\end{equation}
\end{defn}

\begin{rmk}
To see that Definition \ref{def:imm-cent-action} is indeed a rephrasing of Definition \ref{def:imm-cent}, suppose that
  $(G, S, c)$ satisfy the hypotheses of Theorem \ref{gamma cartan exists} (for a trivial 2-cocycle $\omega)$, and consider the associated $H = G/S \rtimes \widehat S, T = \widehat S,$ and actions as in Proposition \ref{prop:weyl-satisfies}. Set  $\Gamma = G/S = H/T$ and $ K =  S = \widehat T$, with the action of $\Gamma$ on $K$ given by Lemma \ref{lem:T-hat-action}.  Then 
  as we showed in the second paragraph of the proof of Theorem \ref{thm:uniqueness},
\[ [\gamma ]  \hact s (\hat t)   = \langle s, [\gamma^{-1}]. \hat t \rangle  = \hat t(\gamma s \gamma^{-1}).\]
In other words, $[\gamma ]  \hact s = \gamma s \gamma^{-1}$.  Comparing Definitions \ref{def:imm-cent-action} and \ref{def:imm-cent}, we see that the action $[\gamma] \hact s$ is immediately centralizing iff $S \leq \Iso(G)$ is immediately centralizing.
\end{rmk}

\begin{thm}
\label{thm:reconstruction-1}
Let $(H, T)$ be as in Assumption \ref{assumption-on-actions}, and assume the existence of $\theta: (H/T) * (H/T) \to \widehat T$ satisfying Equation \eqref{eq:omega-for-recovering-G}. 
Define $ \mathcal G := H/T \rectimes \widehat T$, with the subgroupoid $\mathcal S := \widehat T$; define $\frak s: (\G/ \mathcal S  = H/T) \to \mathcal G$ by $\frak s([\gamma]) = ([\gamma], \id_{p_s(\gamma)})$.  

Assume furthermore  that, for a discrete Abelian group $\Gamma$, we have a homomorphism $\tilde c: H \to \Gamma$ such that $\tilde c$ descends to a homomorphism on $H/T$ and $\tilde c^{-1}(0) $ is effective; and that, for all integers $k \geq 1$, the action of $H/T$ on $\widehat T$ is an immediately centralizing action. 
Then $(\mathcal G, \mathcal S, \frak s)$ satisfies the hypotheses of Theorems \ref{gamma cartan exists} and \ref{thm:Weyl-twist} (with a trivial 2-cocycle $\omega$). In particular,  $C^*_r(\widehat T) \subseteq C^*_r(H/T \rectimes T)$ is $\Gamma$-Cartan.
\end{thm}
\begin{proof}
We first define $c: \G \to \Gamma$ by restricting $\tilde c$:
\[ c([\gamma], \chi) = \tilde c(\gamma),\]
which is well-defined by hypothesis.  Then $c^{-1}(0) = [\tilde c^{-1}(0)] \utimes_\theta \widehat T$.  As $\widehat T$ is a bundle of Abelian groups,  the definition of multiplication in $H/T \rectimes \widehat T$, together  with the fact that $\G$ has the product topology, implies that $\widehat T = \S \subseteq \Iso(c^{-1}(0))^\circ$.  Indeed, our hypothesis that $\tilde c^{-1}(0)$ is effective implies that $\Iso(c^{-1}(0))^\circ = \widehat T$.


To see this, suppose that $U\subseteq \Iso(c^{-1}(0)) $ is an open neighborhood of $([\gamma], \chi)$.  Fix a representative $\gamma$ of $[\gamma]$.  Since the quotient map $H \to H/T$ is continuous, there is an open neighborhood $W$ of $\gamma$ with $[W] = U$.  Therefore, $0 = c(U) =\tilde c(W)$.  As we assumed $c^{-1}(0)$ was effective, we must have $W \subseteq H\z = T$.  Therefore, $([\gamma], \chi) \in \widehat T$.

In particular, as $\S = \widehat T = \Iso(c^{-1}(0))^\circ$ is Abelian, $\S$ is maximal among open Abelian subgroupoids of $\Iso(c^{-1}(0))$.  Moreover, $\S \leq \G$ is closed, since $\G = H/T \utimes_\theta \widehat T$ has the product topology. 

To see that $\S$ is normal, fix $([\gamma], \chi) \in \G, \ ([v], \xi) \in \widehat T$ with $\hat p(\chi)  = \hat p(\xi)$ and $[s(\gamma)] = [v]$, and compute: 
\begin{align*}
([\gamma], \chi) ([v], \xi)([\gamma], \chi)^{-1} &= ([\gamma], \theta([\gamma], [v]) \chi \xi)([\gamma]^{-1}, \theta([\gamma], [\gamma]^{-1})^{-1} ([\gamma]\hact \chi^{-1})) \\ 
 &= ([r(\gamma)], ([\gamma]\hact (\theta([\gamma], [v]) \chi \xi)) \theta ([\gamma], [\gamma]^{-1})^{-1} ([\gamma]\hact \chi^{-1}))  \\
 & = ([r(\gamma)], [\gamma]\hact (\theta([\gamma], [v]) \chi \xi \chi^{-1})\theta ([\gamma], [\gamma]^{-1})^{-1} )\\
 &= ([r(\gamma)], [\gamma]\hact (\theta([\gamma], [v])\xi)\theta ([\gamma], [\gamma]^{-1})^{-1}  ) \in \widehat T
\end{align*}
as desired.

The fact that $\S$ is immediately centralizing comes from the hypothesis that the action of $H/T$ on $\widehat T$ is an immediately centralizing action.   Suppose that there exist $([\gamma], \chi) \in \Iso(\G)$ and an integer $k \geq 1$ such that, for all $\xi \in \widehat T_{p_s(\gamma)} = \widehat T_{p_r(\gamma)}$, there exists $1\leq n \leq k$ with
\begin{equation}
([\gamma], \chi)([s(\gamma)], \xi)^n = ([s(\gamma)], \xi)^n ([\gamma], \chi).
\label{eq:imm-cent-check}
\end{equation}
The fact that $([\gamma], \chi) \in \Iso(\G)$ means that, in particular, $[r(\gamma)] = [s(\gamma)]$.

Since 
units act trivially on $\widehat T$, $([s(\gamma)], \xi)^n = ([s(\gamma)], \theta([s(\gamma)], [s(\gamma)])^{n-1} \xi^n).$  Thus,  Equation \eqref{eq:imm-cent-check} holds iff 
\[([\gamma], \xi^n \chi) = ([\gamma], \theta([s(\gamma)], [\gamma])\theta([s(\gamma)], [s(\gamma)])^{n-1} \chi \xi^n).\]
 The fact that $\widehat T$ is Abelian means that this equation holds, for any $n$ and any $\xi$, iff $\theta([s(\gamma)], [\gamma]) =\theta([r(\gamma)], [\gamma])$ is trivial.
As this last condition is independent of $n$, we conclude that Equation \eqref{eq:imm-cent-check} holds for some $n$ iff it holds for $n=1$.  That is, $\S$ is immediately centralizing.

Finally, the continuity of the section $\frak s$ is an immediate consequence of the fact that $H/T \utimes_\theta \widehat T$ has the product topology. 
\end{proof}

\bibliographystyle{amsalpha}
\bibliography{eagbib}

\end{document}